\documentclass{amsart}
\usepackage{amsmath}
\usepackage{amssymb}
\usepackage{amsthm}
\usepackage{graphicx}
\usepackage{float}
\usepackage{pgfplots}
\usepackage{tikz}
\usetikzlibrary{arrows.meta}
\usepackage{subcaption}
\usepackage{multirow}

\theoremstyle{thmstyleone}%
\newtheorem{theorem}{Theorem}[section]%
\newtheorem{lemma}{Lemma}[section]%

\theoremstyle{thmstyletwo}%
\newtheorem{example}{Example}[section]%
\newtheorem{remark}{Remark}[section]%

\theoremstyle{thmstylethree}%
\newtheorem{definition}{Definition}[section]%

\theoremstyle{thmstylefour}%
\newtheorem{corollary}{Corollary}[section]%

\numberwithin{equation}{section}

\begin{document}

\title[MOSER ITERATION FEASIBILITY FOR GRADIENT ESTIMATES]{The Feasibility of Nash--Moser iteration for Cheng--Yau-type gradient estimates of nonlinear equations on complete Riemannian manifolds}

%    author one information
% \author[short version for running head]{name for top of paper}
\author[B. Shen]{Bin Shen}
\address{School of Mathematics, Southeast University, Nanjing 211189, P. R. China}
%\curraddr{ }
\email{shenbin@seu.edu.cn}
%\thanks{The first author is supported partially by the NNSFC (No. 12001099, 12271093)}

%    author two information
\author[Y. Zhu]{Yuhan Zhu}
\address{School of Mathematics, Southeast University, Nanjing 211189, P. R. China}
%\curraddr{ }
\email{yuhanzhu@seu.edu.cn}

%    \subjclass is required.

%    Abstract is required.
\begin{abstract}
	In this manuscript, we employ the  Nash--Moser iteration technique to determine a condition under which the positive solution $u$ of the generalized nonlinear Poisson equation $$\operatorname{div} (\varphi(|\nabla u|^2)\nabla u) + \psi(u^2)u = 0,$$
	on a complete Riemannian manifold with Ricci curvature bounded from below{,} can be shown to satisfy a Cheng--Yau-type gradient estimate. We define a class of $\varphi$-Laplacian operators by $\Delta_{\varphi}(u):=\operatorname{div} (\varphi(|\nabla u|^2)\nabla u)$, where $\varphi$ is a $C^2$-function under some certain growth conditions. This can be regarded as a natural generalization of the $p$-Laplacian, the $(p,q)$-Laplacian and the exponential Laplacian, as well as having a close connection to the prescribed mean curvature problem. 
	We illustrate the feasibility of applying the Nash--Moser iteration for such Poisson equation to get the Cheng--Yau-type gradient estimates in different cases with various $\varphi$ and $\psi$. Utilizing these estimates, we prove the related Harnack inequalities and a series of Liouville theorems.
	Our results can cover a wide range of quasilinear Laplace operators (\textit{e.g.} $p$-Laplacian for $\varphi(t)=t^{p/2-1}$), and Lichnerowicz-type nonlinear equations (\textit{i.e.} $\psi(t) = At^{p} + Bt^{q} + Ct\log t + D$).
\end{abstract}

\maketitle
%    Text of article.
\section{Introduction}
Let $M$ be an $n$-dimensional complete Riemannian manifold with $\operatorname{Ric}\geq-(n-1)K$ for some $K\geq0$, Cheng and Yau \cite{cheng1975differential} proved  that for a positive harmonic function on geodesic ball $B(o,R)$, there is a constant $c_n$ depending only on $n$ such that
\begin{equation}\label{CY_esti}
\sup_{B(o,R/2)}\frac{|\nabla u|}{u}\leq c_n\frac{1+\sqrt{K}R}{R}.
\end{equation}
This type of gradient estimate is a versatile tool for studying topological and geometrical
properties of manifolds. From (\ref{CY_esti}), for instance, the Harnack inequality, Liouville theorem, estimates of first eigenvalues, as well as optimal Gaussian estimates of the heat kernel can be deduced \cite{li2012geometric}.

One current trend in gradient estimate is to apply {the} Cheng--Yau method to other nonlinear partial differential equations in {the} form of
\begin{align}\label{general_nonlinear_1}
u_t - \Delta u = \Sigma(x,u,t),
\end{align}
with nonlinear function $\Sigma(x,u,t):M\times\mathbb{R}\times[0,+\infty)\to \mathbb{R}$. For example, the classical Li-Yau estimate on {the} Schr\"{o}dinger equation in \cite{li1986parabolic}, the logarithmic type nonlinearities $\Sigma(x,u)= A\log u + Bu$ in \cite{ma2006gradient,yang2008gradient} or the general one in \cite{huang2021gradient,taheri2023soupletzhang}.

Another type of nonlinear equation is
\begin{align}\label{general_nonlinear_2}
\operatorname{div}\left(\mathcal{F}(x,|\nabla u|,\nabla u)\right) = 0,
\end{align}
containing nonlinearity $\mathcal{F}: M \times [0,\infty) \times \Gamma(TM) \to \Gamma(TM)$ inside the divergence operator. Some regularity theorems of such equation in {the} Euclidean space have been investigated by P. Tolksdorf \cite{TOLKSDORF1984126}. In particular, by choosing different $\mathcal{F}$ we got the $p$-Laplacian $\Delta_pu := \operatorname{div}\left(|\nabla u|^{p-2}\nabla u\right)$ and exponential Laplacian $\Delta_eu :=\operatorname{div}\left(\exp\left(1/2|\nabla u|^{2} \right)\nabla u\right)$, of which the Cheng--Yau estimates have been established by B. Kotschwar \cite{kotschwar2009local}, and J. Wu \cite{wu2014gradient} under the condition of a lower bound of sectional curvature. However, when following the traditional Cheng--Yau method for this type of nonlinear equation, the sectional curvature condition becomes necessary due to the use of the Hessian comparison theorem. Then X. Wang \cite{wang2010local} used Nash--Moser iteration technique to weaken the curvature condition, with only a lower bound of Ricci curvature assumed. Moreover, this strategy can be applied not only to more complicated $p$-Laplacian equations (see \cite{han2023gradient,he2023nashmoser,he2023gradient}), but also to more generalized spaces, say, Finsler metric measure space. The Nash--Moser iteration technique is still powerful to bypass the nonlinearity of Finslerian Laplacian $\Delta^{\nabla u}u$, and the analogous Cheng--Yau and Li--Yau estimates have been established by C. Xia \cite{xia2014local} and Q. Xia \cite{xia2022local,xia2023li}.

Nevertheless, not all the gradient estimates obtained in previous research are of the Cheng--Yau type. In contrast, Cheng--Yau estimate is more significant and useful in geometric analysis, since it can derive a strong Liouville property that the bounded positive solution must be a constant (\textit{cf.} Theorem \ref{liou}). So it is of interest to ascertain   \textit{under what circumstances, the differential equation exhibits a Cheng--Yau estimate}.

To answer this question, we shall consider a class of generalized Laplacian operators, motivated by the work of M. Ara \cite{Ara1999GeometryOF}, who introduced the $F$-energy of smooth map $\Phi$ between Riemannian manifolds $(M,g)$ and $(N,h)$ by
\begin{align}
\begin{aligned}
E_F(\Phi):= \int_MF\left(\frac{|d \Phi|^2}{2}\right),
\end{aligned}
\end{align}
where $F:[0,+\infty) \to [0,+\infty)$ is a $C^2$-function with $F'>0$ on $(0,\infty)$. And $F$-harmonic map is defined to be the critical point of $E_F$. For this generalized harmonic map, Y. Dong and his collaborators have already explored vanishing theorem  and Liouville theorem for $F$-harmonic map (or function) \cite{Dong2010OnVT,Dong2015LiouvilleTF}. Now consider the $F$-harmonic function ({\textit{i.e.}} the target manifold is $N = \mathbb{R}$), and the Eular--Langrange equation with respect to $E_F$ is
\begin{align}
\operatorname{div} \left(F'\left(\frac{|\nabla u|^2}{2}\right)\nabla u\right) = 0.
\end{align}
For convenience, we denote such operator by the $\varphi$-Laplacian, namely,
\begin{align}
\Delta_\varphi u:= \operatorname{div} (\varphi(|\nabla u|^2)\nabla u).
\end{align}
Apart from $p$-harmonic or exp-harmonic function, this equation is related to minimal surfaces and prescribed mean curvature \cite{ledesma2023multiplicity}, by setting $\varphi(t) = (1+t)^{-1/2}$, that is,
$$\operatorname{div} \left(\frac{\nabla u}{\sqrt{1+|\nabla u|^2}}\right) = f(u).$$

It is worthwhile to remark that in a weighted Riemannian space, $\varphi$-Laplacian is different from the $f$-Laplacian, which is defined as $\Delta_f u:= e^{f}\operatorname{div} (e^{-f}\nabla u)$
where $f$ is a fixed function independent of the function $u$ to be solved. Also the ``Ricci tensor" in weighted Riemannian space is actually $m$-Bakry-\'{E}mery Ricci tensor  $\operatorname{Ric}^{m,n}_f := \operatorname{Ric} + \operatorname{Hess}(f) - \frac{1}{m-n}df\otimes df$.\\

In the manuscript, we focus on the generalized nonlinear Poisson equation
\begin{align}\label{defharm}
\operatorname{div} (\varphi(|\nabla u|^2)\nabla u) + \psi(u^2)u =0 ,
\end{align}
where $\varphi(t)$ and $\psi(t)$ are $C^\infty$-function on $[0,\infty)$ satisfying $\varphi(t)>0$ for $t>0$.
Equation (\ref{defharm}) arises in the study of reaction-diffusion models with diffusional coefficient $\varphi(|\nabla u|^2)$ and reaction function $\psi(u^2)u$. This equation also has a wide range of applications in physics and engineering.

We present the main theorem as follows.
\begin{theorem}\label{mainthm}
	Let $(M^n,g)$ be a complete Riemannian $n$-manifold with Ricci curvature bounded from below by $\operatorname{Ric} \geqslant -K$ {where} $K\geqslant0$, and let $u$ be a positive solution of (\ref{defharm}) on the ball $B(o,2R)\subset M$. Suppose that for any $t\in [0,\infty)$, $\varphi$ satisfies that
	\begin{equation}\tag{$\varphi_1$}
		-1< l_{\varphi} \leqslant \delta_{\varphi}(t) \leqslant d_{\varphi} < +\infty,
	\end{equation}
	\begin{equation}\tag{$\varphi_2$}
	0< \gamma_{\varphi} \leqslant \frac{\left(\delta_{\varphi}(t)+1\right)^2 }{n-1} - 2t\delta_{\varphi}'(t) \leqslant \varGamma_{\varphi} < +\infty,
	\end{equation}
	where $l_{\varphi}$, $d_{\varphi}$, $\gamma_{\varphi}$, $\varGamma_{\varphi}$ are all constants, and  
	\begin{eqnarray}\label{deltavarphi}
	\delta_{\varphi}(t) := \frac{2t\varphi'(t)}{\varphi(t)}.
	\end{eqnarray}
	Moreover, $\delta_{\psi}(t)$ is defined in the same way by replacing $\varphi$ by $\psi$ in \eqref{deltavarphi}, satisfying \begin{align}\label{psi2}
	\varTheta_{\varphi,\psi} := \sup_{\substack{s\geqslant 0,\\t\in \mathbb{R}^+ -  I_{\psi}}} \left(\frac{2\left(\delta_{\varphi}(s)+1\right) }{n-1} +\delta_{\varphi}(s) - \delta_{\psi}(t) \right)^2  < \frac{4\gamma_{\varphi}}{n-1},
	\end{align}
	where
	\begin{align}\label{psi1}
	I_{\psi}:= \left\{t>0: \psi(t)\left[\frac{2\left(\delta_\varphi(s) +1 \right)}{n-1} + \delta_\varphi(s)-\delta_\psi(t) \right] \geqslant 0, \text{ for each }s\geqslant 0\right\}.
	\end{align}
	%    \begin{itemize}
	%    	\item[$(\psi_1)$] 
	%    	$\psi(t) \leqslant 0,~~ \delta_{\psi}(t)\geqslant\frac{2(d_{\varphi}+1)}{n-1}+d_{\varphi},$
	%    	\item [$(\psi_2)$] 
	%    	$\psi(t) \geqslant 0,~~ \delta_{\psi}(t)\leqslant\frac{2(l_{\varphi}+1)}{n-1}+l_{\varphi},$
	%    	\item [$(\psi_3)$] 
	%    	$\sup_{t,s\geqslant 0}\left(\frac{2\left(\delta_{\varphi}(t)+1\right) }{n-1} +\delta_{\varphi}(t) - \delta_{\psi}(s) \right)^2  < \frac{4\gamma_{\varphi}}{n-1}.$
	%    \end{itemize}
	Then, there exists a constant $C= C(n,l_{\varphi},d_{\varphi},\gamma_{\varphi},\varGamma_{\varphi},\varTheta_{\varphi,\psi})$ which depends only on $n$ and those constants related to the equation itself (in fact, $\delta_\varphi$ and $\delta_\psi$), such that
	\begin{align}\label{esti}
	\frac{|\nabla u|}{u}\leqslant C\frac{1+\sqrt{K}R}{R}
	\end{align}
	on $B(o,R)$.
\end{theorem}

\begin{remark}
	When $I_\psi = (0,+\infty)$, condition (\ref{psi2}) is naturally satisfied since $\mathbb{R}^+ - I_\psi = \emptyset$. Consequently, the constant $C$ in the estimate (\ref{esti}) depends only on the diffusional coefficient $\varphi$.
\end{remark}

\begin{remark}
	In the case that $u$ is negative, one may consider $-u$ as a positive solution of euqation (\ref{defharm}). Hence, the estimate (\ref{esti}) is also valid for negative solutions.
\end{remark}

Based on the aforementioned estimate, we show some immediate consequences of Theorem \ref{mainthm} as follows.
\begin{theorem}[Harnack's inequality]\label{harnack}
	Under the same assumption in Theorem \ref{mainthm}, there exists a constant $C(n,l_{\varphi},d_{\varphi},\gamma_{\varphi},\varGamma_{\varphi},\varTheta_{\varphi,\psi})$ such that for any $x, y \in B(R)$,
	$$
	u(x) / u(y) \leq e^{C(1+\sqrt{K} R)} .
	$$
	
	It follows that if $K=0$, then we have a constant independent of $R$ such that
	$$
	\sup _{B(R)} u \leq C\inf _{B(R)} u .
	$$
	
\end{theorem}
\begin{theorem}[Liouville theorem]\label{liou}
	Let $M$ be a complete and non-compact Riemannian manifold with non-negative Ricci curvature, and let $u$ be a bounded positive solution of (\ref{defharm}) with $\varphi$ and $\psi$ satisfying the same assumption in Theorem \ref{mainthm}. If $\psi(t) = 0$ has positive root $t=T$ then $u \equiv \sqrt{T}$. Otherwise, there is no such positive solution.
\end{theorem}

To interpret conditions ($\varphi_1$) and ($\varphi_2$), we take the $(p,q)$-Laplacian as a non-trivial example, which generalizes the results for $p$-Laplacian.
\begin{corollary}\label{pq}
	Let $u$ be a positive solution of the following equations
	$$\Delta_{p,q}u := \operatorname{div}\left(\left(|\nabla u|^{p-2} + |\nabla u|^{q-2} \right)\nabla u  \right) = 0$$
	on $B(o,2R)$. If $p,q>1$ and 
	\begin{align}\label{req_pq}
	(n-1)<\frac{4(p-1)(q-1)}{(p-q)^2},
	\end{align}
	then we have
	\begin{align*}
	\frac{|\nabla u|}{u}\leqslant C(n,p,q)\frac{1+\sqrt{K}R}{R}
	\end{align*}
	on $B(o,R)$.
\end{corollary}

Or more generally, we can consider a finite linear combination of several $p_i$-Laplacian operators, called weighted $(p_1,...,p_r)$-Laplacian.
\begin{corollary}\label{pi}
	Let $u$ be a positive solution of
	$$\tilde{\Delta}_{p_1,...,p_r}u:=\left(\sum_{i=1}^ra_i\Delta_{p_i}\right)u = \operatorname{div}\left(\sum_{i=1}^{r}a_i|\nabla u|^{p_i-2}\nabla u  \right) = 0$$
	on  $B(o,2R)$, where $a_i >0$ and $1<p_1<p_2<\cdots<p_r$.
	If 
	\begin{align}\label{req_pi}
	(n-1)<\frac{2(p_1-1)^2}{(p_r-p_1)^2},
	\end{align}
	then we have
	\begin{align*}
	\frac{|\nabla u|}{u}\leqslant C(n,r,p_i)\frac{1+\sqrt{K}R}{R}
	\end{align*}
	on $B(o,R)$.
\end{corollary}
We will give a detailed explanation in Section 4. Here, to sum up, we list some common and new Laplacian operators in the following table.
\begin{table}[H]
	\centering
	\caption{Related constants for different Laplacian operators}
	\begin{tabular}{|c|c|c|c|c|}
		\hline
		& $\Delta$ & $\Delta_p$ & $\Delta_{p,q}$ &\begin{tabular}[c]{@{}c@{}} $\tilde{\Delta}_{p_1,...,p_r}$ \\($p_1<...<p_r$) \end{tabular}\\\hline
		$\varphi(t)$    &1      &$t^{p/2-1}$         &$t^{p/2-1}+t^{q/2-1} $  &$\sum_{i=1}^{r}a_it^{p_i/2-1}$                                \\\hline
		$\delta_\varphi(t)$    &0          & $p-2$           &$\frac{(p-2)t^{p/2-1}+(q-2)t^{q/2-1} }{t^{p/2-1}+t^{q/2-1}}$                  &$\frac{\sum_{i=1}^{r}a_i(p_i-2)t^{p_i/2-1}}{\sum_{i=1}^{r}a_it^{p_i/2-1}}$                                \\\hline
		$d_\varphi$         &0          & $p-2$             &$\max\{p,q\}-2$                &${p_r}-2$                                \\\hline
		$l_\varphi$         &0          & $p-2$             &$\min\{p,q\}-2$                &${p_1}-2$                                \\\hline
		$\gamma_\varphi$    &$\frac{1}{n-1}$          &$\frac{(p-1)^2}{n-1}$            &$\frac{4(p-1)(q-1)-(n-1)(q-p)^2}{4n}$                &$\frac{({p_1}-1)^2}{n-1}-\frac{(p_r-p_1)^2}{2}$                                \\\hline
		$\varGamma_\varphi$ &$\frac{1}{n-1}$           & $\frac{(p-1)^2}{n-1}$           &$\frac{(\max\{p,q\}-1)^2}{n-1}$                &$\frac{(p_r-1)^2}{n-1}$                               \\\hline
	\end{tabular}
\end{table}
\begin{remark}
	When $r=2$, Corollary \ref{pi} reduces to Corollary \ref{pq}, {whereas} the condition {(\ref{req_pi}) for $p_i$} will be slightly stronger than {(\ref{req_pq})} in Corollary \ref{pq}. In fact, (\ref{req_pq}) is the sufficient and necessary condition for the existence of positive $\gamma_{\varphi}$. Due to the lack of explicit solution for the high-degree polynomial  equations, it is unlikely to find a precise infimum as what we did in Corollary \ref{pq} (\textit{cf.} Example \ref{ex_pq} and Example \ref{ex_pr}).
\end{remark}

Next we take some special cases of $\psi$ and $\varphi$ in equation (\ref{defharm}){,} in order to compare our results with those obtained in previous research.
When $\varphi\equiv1$ and $\psi(t)=1-t$, (\ref{defharm}) becomes Allen--Cahn equation
$$\Delta u + (1-u^2)u =0.$$ Theorem \ref{mainthm} {improves} the result in \cite{hou2019gradient}{,} as we do not need the bounded condition $u\leqslant C${, so that} the estimate is independent of the upper bound of $u$. Also, our result is exactly Cheng--Yau estimate, without the correction term $(1-u^2)$. More generally, when $\psi(t)=t^{(m-1)/2}-t^{(k-1)/2}$, that is,
$$\Delta u + u^{m}-u^k =0,$$
Y. Wang \cite{wang2024gradient} has proved a Liouville Theorem for
$$1<m<\frac{n+3}{n-1} \text{ or } 1<k<\frac{n+3}{n-1},$$ whereas our outcome has weaker conditions (see Remark \ref{pic_umus}) and can generalize \cite{wang2024gradient} to $p$-Laplacian or even $(p_1,...,p_r)$-Laplacian (see Table 2).
Recently, J. He and Y. Wang \cite{he2023nashmoser} also studied the generalized Lane--Emden equation 
$$\Delta_{p}u + au^q = 0,$$
which means $\varphi(t) = t^{p/2-1}$ and $\psi(t)=at^{(q-1)/2}$ in (\ref{defharm}). This equation is also related to {prescribed} scalar curvature problem that 
$$\Delta u + Ku^{\frac{n+2}{n-2}} = 0.$$
Theorem \ref{mainthm} shows the Cheng--Yau estimate when 
$$q<\frac{n+3}{n-1}(p-1), a>0,$$ or 
$$q>p-1, a<0,$$
covering the result obtained in \cite{he2023nashmoser}. Moreover, our result can be extended to more general equations. We list the brief results in the table as follows (see details in Section 4).
\begin{table}[H]
	\centering
	\caption{Liouville Theorems for different Laplacian equations $\Delta_\varphi + \psi(u^2)u$ = 0}
	\begin{tabular}{|c|c|c|}
		\hline
		& $au^q$ & $u^m-u^k ~~(m<k)$   \\ \hline
		$\Delta$  & \begin{tabular}[c]{@{}c@{}}$a>0 \text{ and }{q}<\frac{n+3}{n-1}$\\ or\\ $a<0 \text{ and }{q}>1$ \\\\ ($\Rightarrow$ No positive bounded solution)\end{tabular}   & \begin{tabular}[c]{@{}c@{}}$m<\frac{n+3}{n-1}$\\ and\\ $k>1$ \\\\ ($\Rightarrow$ $u\equiv1$)\end{tabular}           \\ \hline
		$\Delta_p$ & \begin{tabular}[c]{@{}c@{}}$a>0 \text{ and }\frac{q}{p-1}<\frac{n+3}{n-1}$\\ or\\ $a<0 \text{ and }{q}>{p-1}$\\\\ ($\Rightarrow$ No positive bounded solution)\end{tabular} & \begin{tabular}[c]{@{}c@{}}$\frac{m}{p-1}<\frac{n+3}{n-1}$\\ and\\ $\frac{k}{p-1}>1$\\\\ ($\Rightarrow$ $u\equiv1$)\end{tabular}          \\ \hline
		\begin{tabular}[c]{@{}c@{}}$\tilde{\Delta}_{p_1,...,p_r}$\\{$n<\mathcal{N}_1$}\end{tabular} & \multirow{2}{*}{\begin{tabular}[c]{@{}c@{}}$a>0$ and\\$ \frac{q-\gamma}{p_1-1}<\frac{n+1}{n-1}$\\ or\\ $a<0$ and\\$\frac{q+\gamma}{p_r-1}>\frac{n+1}{n-1}$\\\\where \\$\gamma = 2\sqrt{\frac{(p_1-1)^2}{(n-1)^2}-\frac{(p_r-p_1)^2}{2(n-1)}}$\\\\ ($\Rightarrow$ No positive bounded solution)\end{tabular}} & \begin{tabular}[c]{@{}c@{}}{$\frac{m}{p_1-1}\leqslant\frac{n+1}{n-1}$}\\ {and}\\ {$\frac{k}{p_r-1}\geqslant\frac{n+1}{n-1}$}\\\\ ($\Rightarrow$ $u\equiv1$)\end{tabular}           \\ \cline{1-1} \cline{3-3}
		\begin{tabular}[c]{@{}c@{}}$\tilde{\Delta}_{p_1,...,p_r}$\\{$n<\mathcal{N}_2$}\end{tabular} &  & \begin{tabular}[c]{@{}c@{}}$\frac{m-\gamma}{p_1-1}<\frac{n+1}{n-1}$\\ and\\ $\frac{k+\gamma}{p_r-1}>\frac{n+1}{n-1}$\\\\where \\$\gamma = 2\sqrt{\frac{(p_1-1)^2}{(n-1)^2}-\frac{(p_r-p_1)^2}{2(n-1)}}$\\\\ ($\Rightarrow$ $u\equiv1$)\end{tabular}           \\ \hline
	\end{tabular}
\end{table}

{In the last two rows of Table 2, $\mathcal{N}_1(p_1,p_r) > \mathcal{N}_2(p_1,p_r)$ are the first and second critical dimensions of $\tilde{\Delta}_{p_1,...,p_r}$, defined by
$$
\mathcal{N}_{1}:= 2\left(\frac{\min\{p_i\} - 1}{\max\{p_i\}- \min\{p_i\}}\right)^2 +1,
$$
and
$$
\mathcal{N}_{2}:= \sqrt{2\mathcal{N}_1 +3} -2,
$$
respectively, depending only on the maximum and minimum of $p_i$. The dimension $n$ with different critical conditions 
can lead to different Liouville properties as shown in that table.}

Furthermore, $\psi$ is not necessarily a polynomial, especially, we now consider \begin{align*}
\Delta u + uh(\log u) = 0. 
\end{align*}Although there has been several research on this equation 
(\textit{cf.} \cite{ma2006gradient,huang2021gradient,peng2021yau}), these estimates are not of the Cheng--Yau type, and therefore cannot be used to prove Liouville theorems like Theorem \ref{liou}. Here we present our new results in table below (refer to Corollary \ref{log} for details). 
\begin{table}[H]
	\centering
	\caption{Liouville Theorems for different Laplacian equations with logarithm}
	\begin{tabular}{|c|c|}
		\hline
		& \begin{tabular}[c]{@{}c@{}}$\Delta_\varphi + au^q(\log u)^m = 0$ \\\\($m$ is rational with \\$m=(2k_1+1)/(2k_2+1)$ and $am<0$)\end{tabular}    \\ \hline
		$\Delta$  & \begin{tabular}[c]{@{}c@{}}$1<q<\frac{n+3}{n-1}$\\\\ ($\Rightarrow$ $u\equiv 1$ when $m>0$ \\No positive bounded solution when $m<0$)\end{tabular}              \\ \hline
		$\Delta_p$ & \begin{tabular}[c]{@{}c@{}}${p-1}<\frac{q}{p-1}<\frac{n+3}{n-1}$\\ \\ ($\Rightarrow$ $u\equiv 1$ when $m>0$ \\No positive bounded solution when $m<0$)\end{tabular}          \\ \hline
		\begin{tabular}[c]{@{}c@{}}$\tilde{\Delta}_{p_1,...,p_r}$\\ $(a_1>0)$\end{tabular} & \begin{tabular}[c]{@{}c@{}}$ \frac{q}{p_1-1}<\frac{n+1}{n-1}+2\sqrt{\frac{1}{(n-1)^2}-\frac{(p_r-p_1)^2}{2(n-1)(p_1-1)^2}}$ \\and\\$\frac{q}{p_r-1}>\frac{n+1}{n-1}-2\sqrt{\frac{(p_1-1)^2}{(n-1)^2(p_r-1)^2}-\frac{(p_r-p_1)^2}{2(n-1)(p_r-1)^2}}$\\and\\ $\left(\frac{(n+1)^2}{4} + \frac{n-1}{2} \right)(p_r-p_1)<(p_1-1)^2$\\\\ ($\Rightarrow$ $u\equiv 1$ when $m>0$ \\No positive bounded solution when $m<0$)\end{tabular}            \\ \hline
	\end{tabular}
\end{table}

This manuscript is arranged as follows. In Section 2, we introduce basic definitions, {and} derive a Bochner-type formula which is a necessary tool in {the} Moser's iteration. Then we prove the main theorem in Section 3. In {the} last section, we discuss some specific examples.

\section{Preliminary}

We consider a (weak) positive  solution $u\in C^1(\Omega)\cap W^{1,p}(\Omega)$ of equation (\ref{defharm}) over a bounded reign $\Omega\subset M$, which means
\begin{align}
- \int_{\Omega}\varphi(|\nabla u|^2)\langle \nabla u, \nabla \phi \rangle + \int_{\Omega}\psi(u^2)u\phi = 0 
\end{align}
for any $v\in C^\infty_0(\Omega)$.

Let $M_\varepsilon := \{x\in M: |\nabla u| (x) > \varepsilon/2\}$ for some $\varepsilon >0$. Since $\varphi(t)=0$ only holds for $t=0$, 
$$\inf_{t\in [\frac{\varepsilon^2}{4}, C]}\varphi(t)>0$$
for any fixed $C>0$. The regularity theorem (see Remark 2.7 in \cite{cianchi2018second}) shows that $u\in W^{2,2}_{\text{loc}}(\Omega\cap M_\varepsilon)$. It should be noted that $\Delta_{\varphi}$ is not necessarilly a uniformly elliptic operator in $M$. However, in the regular part $M_\varepsilon$.  $\Delta_{\varphi}$ is uniformly elliptic. Hence, by usual bootstrap argument, the weak  solution $u$ is in fact smooth in $\Omega\cap M_\varepsilon$ since both $\varphi(t)$ and $\psi(t)$ are smooth.

Denote $H := |\nabla u|^2$, then (\ref{defharm}) reduces to
\begin{align}\label{defharm2}
\begin{aligned}
\Delta_{\varphi}u &= \operatorname{div} (\varphi(H)\nabla u)\\
&= \varphi(H)\Delta u + \varphi'(H)\left\langle \nabla H, \nabla u \right\rangle \\
&= -\psi(u^2)u.
\end{aligned}
\end{align}
Note that $\varphi$-Laplacian $\Delta_{\varphi}$ is not necessarily a linear operator,  so we shall choose a suitable linearization operator  $\mathcal{L}_\varphi$ defined by
\begin{align}\label{linear}
\begin{aligned}
\mathcal{L}_\varphi(\eta) :=& \operatorname{div}\left({\varphi(H)}\nabla\eta + {2\varphi'(H)} \left\langle \nabla u, \nabla \eta \right\rangle \nabla u \right) \\
=& \operatorname{div}\left(\varphi(H)\mathcal{A}(\nabla \eta) \right),
\end{aligned}
\end{align}
where
\begin{align*}
\mathcal{A} := \operatorname{id} + \frac{2\varphi'(H)\nabla u \otimes \nabla u}{\varphi(H)}.
\end{align*}

Suppose that $u$ is positive and set $f := \log u$ and $\hat{H} := H/u^2 = |\nabla f|^2$. It is easy to check that
\begin{align}\label{u2f_nabla}
\nabla H = {u^2}\left(\nabla \hat{H} + 2\hat{H}\nabla f \right),
\end{align}
and
\begin{align}\label{u2f_Laplacian}
\Delta u = {u}\left(\Delta f + \hat{H} \right).
\end{align}
Hence, (\ref{defharm2}) can be written as
\begin{align}\label{u2f_harmonic}
\begin{aligned}
\varphi(H)\Delta f &= \varphi(H)\left(\frac{\Delta u}{u} - \frac{|\nabla u^2|}{u^2} \right)\\
&= - {\varphi'(H)}\left\langle\nabla H, \nabla f\right\rangle - {\psi(u^2)} -\varphi(H)\hat{H}\\
&= - {\varphi'(H)H}\frac{\left\langle\nabla \hat{H}, \nabla f\right\rangle}{\hat{H}} - \left({2\varphi'(H)H} + {\varphi(H)}\right)\hat{H}  - {\psi(u^2)}.
\end{aligned}
\end{align}

We now define a special function to {simplify (\ref{linear}) and (\ref{u2f_harmonic})}.
\begin{definition}\label{degree}
	For a $C^1$-function $\varphi$ on $[0,\infty)$, {the} \textbf{degree function} {of $\Delta_\varphi$} is defined by
	$$\delta_\varphi(t) := \frac{2t\varphi'(t)}{\varphi(t)}.$$
	We say $\varphi$ has finite lower degree and upper degree if there exist finite constants such that
	$$
	\inf_{t\geqslant 0}\delta_\varphi(t) = l_{\varphi} > -\infty,~~~~\sup_{t\geqslant 0}\delta_\varphi(t) = d_{\varphi} < +\infty.
	$$
\end{definition}

Then we will derive the following Bochner-type formula for the linearization operator of $\varphi$-Laplacian.
\begin{lemma}\label{btf_l}
	By adopting the same notations as above, any positive $\varphi$-harmonic function $u$ satisfies that 
	\begin{align*}
	\begin{aligned}
	\mathcal{L}_\varphi(\hat{H}) &= 2\varphi(H)\left( |\nabla^2 f|^2  + \operatorname{Ric}(\nabla f)\right) + \varphi'(H)u^2|\nabla  \hat{H}|^2 - 4\delta_\varphi'(H)H\varphi(H) \hat{H}^2\\
	&~~~~+ 2\psi(u^2)\left(\delta_\varphi(H)-\delta_\psi(u^2)\right)\hat{H} - 2\varphi(H)\left(\delta_\varphi(H)+1+\delta_\varphi'(H)H\right)\left\langle \nabla \hat{H}, \nabla f\right\rangle
	\end{aligned}
	\end{align*}
	on $M_\varepsilon$, where $\nabla^2$ denotes the Hessian operator.
\end{lemma}

\begin{proof}
	By the definition of linearization in (\ref{linear}),
	\begin{align}\label{btf_e1}
	\begin{aligned}
	\mathcal{L}_\varphi(\hat{H}) &= \operatorname{div}\left(\varphi(H)\nabla \hat{H} + 2\varphi'(H) \left\langle \nabla u, \nabla \hat{H} \right\rangle \nabla u \right) \\
	&= \operatorname{div}\left(\varphi(H)\nabla \hat{H} + 2\varphi'(H)H\frac{\left\langle \nabla f, \nabla \hat{H} \right\rangle }{\hat{H}}  \nabla f \right) \\
	&= \varphi(H)\Delta \hat{H} + \varphi'(H)\left\langle\nabla H, \nabla \hat{H} \right\rangle + 2\varphi'(H)H\frac{\left\langle \nabla f, \nabla \hat{H} \right\rangle }{\hat{H}}\Delta f  \\
	&~~~~~~+2\left\langle \nabla\left( \varphi'(H)H\frac{\left\langle \nabla f, \nabla \hat{H} \right\rangle }{\hat{H}} \right), \nabla f\right\rangle.
	\end{aligned}
	\end{align}
	Utilizing (\ref{u2f_nabla}), (\ref{u2f_harmonic}) and the standard Bochner formula of Laplacian, namely,
	\begin{align*}
	\frac{1}{2}\Delta \hat{H} = |\nabla^2 f|^2 + \left\langle \nabla\Delta f, \nabla f\right\rangle + \operatorname{Ric}(\nabla f),
	\end{align*}
	one may find that
	\begin{align}\label{btf_e2}
	\begin{aligned}
	\mathcal{L}_\varphi(\hat{H}) =& 2\varphi(H)\left( |\nabla^2 f|^2 + \left\langle \nabla\Delta f, \nabla f\right\rangle + \operatorname{Ric}(\nabla f)\right)+ \varphi'(H)u^2|\nabla  \hat{H}|^2 \\
	& + 2\varphi'(H) {H}\left\langle\nabla f, \nabla \hat{H} \right\rangle+ 2\varphi'(H)H\frac{\left\langle \nabla f, \nabla \hat{H} \right\rangle }{\hat{H}}\Delta f\\
	&+ 2\left\langle \nabla\left( \varphi'(H)H\frac{\left\langle \nabla f, \nabla \hat{H} \right\rangle }{\hat{H}} \right), \nabla f\right\rangle\\
	=& 2\varphi(H)\left( |\nabla^2 f|^2 + \left\langle \nabla\Delta f, \nabla f\right\rangle + \operatorname{Ric}(\nabla f)\right) + \varphi'(H)u^2|\nabla  \hat{H}|^2\\
	& + 2\varphi'(H) {H}\left\langle\nabla f, \nabla \hat{H} \right\rangle+ 2\varphi'(H)H\frac{\left\langle \nabla f, \nabla \hat{H} \right\rangle }{\hat{H}}\Delta f   \\
	&- 2\left\langle \nabla\varphi(H)\Delta f, \nabla f\right\rangle- 2\left\langle \nabla\left(\left(\varphi(H) + 2\varphi'(H)H\right)\hat{H}\right), \nabla f\right\rangle\\
	&- 2\left\langle \nabla\psi(u^2), \nabla f\right\rangle.
	\end{aligned}
	\end{align}
	We can calculate directly the last three terms on the RHS of (\ref{btf_e2}) as the follows.
	\begin{align}\label{btf_e21}
	\begin{aligned}
	&- 2\langle \nabla\varphi(H)\Delta f, \nabla f\rangle\\=& -2\varphi(H)\left\langle\nabla\Delta f , \nabla f\right\rangle- 2\varphi'(H)\Delta f\left\langle\nabla H, \nabla f\right\rangle\\
	=& -2\varphi(H)\left\langle\nabla\Delta f , \nabla f\right\rangle- 2\varphi'(H)H\Delta f\frac{\left\langle \nabla \hat{H}, \nabla f  \right\rangle }{\hat{H}} - 4\varphi'(H)H\Delta f \hat{H}\\
	=& -2\varphi(H)\left\langle\nabla\Delta f , \nabla f\right\rangle- 2\varphi'(H)H\Delta f\frac{\left\langle \nabla \hat{H}, \nabla f  \right\rangle }{\hat{H}}\\
	&+ 4\varphi'(H)H\left( \frac{\varphi'(H)H}{\varphi(H)}{\left\langle\nabla \hat{H}, \nabla f\right\rangle} + \left(\frac{2\varphi'(H)H}{\varphi(H)} + 1\right)\hat{H}^2 + \frac{\psi(u^2)}{\varphi(H)}\hat{H} \right),
	\end{aligned}
	\end{align}
	\begin{align}\label{btf_e22}
	\begin{aligned}
	&- 2\left\langle \nabla\left(\left(\varphi(H) + 2\varphi'(H)H\right)\hat{H}\right), \nabla f\right\rangle\\
	=& -2\left(\varphi(H) + 2\varphi'(H)H\right)\langle \nabla\hat{H}, \nabla f\rangle -2\left(\varphi(H) + 2\varphi'(H)H\right)'\hat{H}\left\langle \nabla H, \nabla f\right\rangle\\
	=& -2\left(\varphi(H) + 2\varphi'(H)H\right)\langle \nabla\hat{H}, \nabla f\rangle -2\left(\varphi(H) + 2\varphi'(H)H\right)'H\langle \nabla \hat{H}, \nabla f\rangle\\
	& - 4\left(\varphi(H) + 2\varphi'(H)H\right)'H\hat{H}^2,
	\end{aligned}
	\end{align}
	and
	\begin{align}\label{btf_e23}
	\begin{aligned}
	- 2\left\langle \nabla\psi(u^2), \nabla f\right\rangle = -4\psi'(u^2)u\left\langle \nabla u, \nabla f\right\rangle = -4\psi'(u^2)u^2\hat{H}.
	\end{aligned}
	\end{align}
	Thus, it follows from (\ref{btf_e2}) -- (\ref{btf_e22}) that
	\begin{align}\label{btf_e3}
	\begin{aligned}
	\mathcal{L}_\varphi(\hat{H})= & 2\varphi(H)\left( |\nabla^2 f|^2  + \operatorname{Ric}(\nabla f)\right)\\
	& + \varphi'(H)u^2|\nabla  \hat{H}|^2+ 4\left(\frac{\varphi'(H)H}{\varphi(H)}\psi(u^2)-\psi'(u^2)u^2 \right)\hat{H}\\
	&+ 2\left[4\left(\frac{\varphi'(H)^2H^2}{\varphi(H)}\right) + 2\varphi'(H)H - 2\left(\varphi(H) + \varphi'(H)H\right)'H\right] \hat{H}^2\\
	&+ 2\left[\frac{2\varphi'(H)^2H^2}{\varphi(H)}-\varphi(H) - \varphi'(H)H - \left(\varphi(H) + \varphi'(H)H\right)'H\right]\langle \nabla \hat{H}, \nabla f\rangle.
	\end{aligned}
	\end{align}
	In term of the degree functions $\delta_\varphi$ and $\delta_\psi$ in Definition \ref{degree},
	we see
	\begin{align*}
	\begin{aligned}
	\left(\varphi(H) + 2\varphi'(H)H\right)'H &= \left(\left(\delta_\varphi(H)+1 \right)\varphi(H)\right)'H \\
	&= \delta'_\varphi(H)H\varphi(H) + \frac{1}{2}\delta_\varphi(H)^2\varphi(H) + \frac{1}{2}\delta_\varphi(H)\varphi(H).
	\end{aligned}
	\end{align*}
	Hence (\ref{btf_e3}) becomes
	\begin{align}\label{btf_final}
	\begin{aligned}
	\mathcal{L}_\varphi(\hat{H}) &= 2\varphi(H)\left( |\nabla^2 f|^2  + \operatorname{Ric}(\nabla f)\right) + \varphi'(H)u^2|\nabla  \hat{H}|^2 \\
	&~~~~+ 2\psi(u^2)\left(\delta_\varphi(H)-\delta_\psi(u^2)\right)\hat{H}- 4\delta_\varphi'(H)H\varphi(H) \hat{H}^2\\
	&~~~~ - 2\varphi(H)\left(\delta_\varphi(H)+1+\delta_\varphi'(H)H\right)\left\langle \nabla \hat{H}, \nabla f\right\rangle.
	\end{aligned}
	\end{align}
	It finishes the proof.
\end{proof}

With the assistance of the Bochner-type formula in Lemma \ref{btf_l}, a lower estimate of $\mathcal{L}_\varphi(\hat{H})$ can be derived as follows.
\begin{lemma}\label{LH_lower}
	A lower bound of $\mathcal{L}_\varphi(\hat{H})$ can be given by
	\begin{align}\label{eq-LH_lower}
	\begin{aligned}
	\mathcal{L}_\varphi(\hat{H}) &\geqslant 2\varphi(H)\operatorname{Ric}(\nabla f) + \varphi(H)(\delta_\varphi(H)+1)\frac{|\nabla \hat{H}|^2}{\hat{H}}\\
	&~~~~ + \left[\frac{2(\delta_\varphi(H)+1)^2}{n-1} - 4\delta_\varphi'(H)H\right]\varphi(H) \hat{H}^2\\
	&~~~~ + \left[\frac{2(\delta_\varphi(H)+1)^2}{n-1}-2(\delta_\varphi(H)+1) -2\delta_\varphi'(H)H\right]\varphi(H)\left\langle \nabla \hat{H}, \nabla f\right\rangle\\
	&~~~~ + \frac{2\varphi(H)}{n-1}\left((\delta_{\varphi}(H)+1)\frac{\left\langle\nabla \hat{H}, \nabla f\right\rangle}{2\hat{H}} + \frac{\psi(u^2)}{\varphi(H)} \right)^2 \\
	&~~~~ + 2\psi(u^2)\left[\frac{2\left(\delta_\varphi(H) +1 \right)}{n-1} + \delta_\varphi(H)-\delta_\psi(u^2) \right]\hat{H}.
	\end{aligned}
	\end{align}
\end{lemma}
\begin{proof}
	We need to estimate the Hessian term $|\nabla^2f|^2$ subtly. Choose a local orthonormal frame $\{e_i\}$ with $e_1 = \nabla f/|\nabla f|$. Then
	\begin{align}\label{f11}
	f_{11} = \frac{\left\langle\nabla \hat{H}, \nabla f\right\rangle}{2\hat{H}}
	\end{align}
	and
	\begin{align}\label{sum_f1i}
	\sum_{i=1}^{n}f_{1i}^2 = \frac{|\nabla \hat{H}|^2}{4\hat{H}}.
	\end{align}
	In such an orthonormal frame, one could immediately deduce from (\ref{u2f_harmonic}) that
	\begin{align}\label{u2f_harmonic_orth}
	\begin{aligned}
	\sum_{i=2}^{n}f_{ii} &= -f_{11} -  \left(\frac{2\varphi'(H)H}{\varphi(H)}\right)f_{11} - \left(\frac{2\varphi'(H)H}{\varphi(H)}+1\right)\hat{H} - \frac{\psi(u^2)}{\varphi(H)}\\
	&=-\left(\frac{2\varphi'(H)H}{\varphi(H)}+1\right)\left(f_{11}+\hat{H}\right)- \frac{\psi(u^2)}{\varphi(H)}\\
	&=-\left(\delta_{\varphi}(H)+1\right)f_{11} - \left(\delta_{\varphi}(H)+1\right)\hat{H}- \frac{\psi(u^2)}{\varphi(H)}.
	\end{aligned}
	\end{align}
	Therefore,
	\begin{align*}
	\begin{aligned}
	|\nabla^2 f|^2 &\geqslant \sum_{i=1}^{n}f_{1i}^2 + \sum_{i=2}^{n}f_{ii}^2 \\
	&\geqslant  \sum_{i=1}^{n}f_{1i}^2 + \frac{1}{n-1}\left(\sum_{i=2}^{n}f_{ii} \right)^2\\
	&\geqslant  \sum_{i=1}^{n}f_{1i}^2 + \frac{1}{n-1}\left(\left(\delta_{\varphi}(H)+1\right)f_{11} + \left(\delta_{\varphi}(H)+1\right)\hat{H}+ \frac{\psi(u^2)}{\varphi(H)}\right)^2.
	\end{aligned}
	\end{align*}
	Since (\ref{f11}) and (\ref{sum_f1i}), it infers that
	\begin{align}\label{esti_hess}
	\begin{aligned}
	|\nabla^2 f|^2 \geqslant \frac{|\nabla \hat{H}|^2}{4\hat{H}} &+ \frac{\left(\delta_\varphi(H) +1 \right)^2}{n-1}\left\langle\nabla \hat{H}, \nabla f\right\rangle + \frac{\left(\delta_\varphi(H) +1 \right)^2}{n-1}\hat{H}^2\\
	&+ \frac{2\left(\delta_\varphi(H) +1 \right)\psi(u^2)}{(n-1)\varphi(H)}\hat{H}+\frac{\left(\delta_\varphi(H) +1 \right)^2}{n-1}f_{11}^2\\
	& + \frac{2\left(\delta_\varphi(H) +1 \right)f_{11}\psi(u^2)}{(n-1)\varphi(H)} + \frac{\psi(u^2)^2}{\varphi(H)^2}.
	\end{aligned}
	\end{align}
	
	Substituting the estimate (\ref{esti_hess}) for $|\nabla^2 f|^2$ in (\ref{btf_final}) yields
	\begin{align*}
	\begin{aligned}
	\mathcal{L}_\varphi(\hat{H}) &\geqslant 2\varphi(H)\operatorname{Ric}(\nabla f) + \varphi(H)(\delta_\varphi(H)+1)\frac{|\nabla \hat{H}|^2}{2\hat{H}} \\
	&~~~~ + \left[\frac{2(\delta_\varphi(H)+1)^2}{n-1}-2(\delta_\varphi(H)+1) -2\delta_\varphi'(H)H\right]\varphi(H)\left\langle \nabla \hat{H}, \nabla f\right\rangle\\
	&~~~~ + \left[\frac{2(\delta_\varphi(H)+1)^2}{n-1} - 4\delta_\varphi'(H)H\right]\varphi(H) \hat{H}^2\\
	&~~~~ + 2\psi(u^2)\left[\frac{2\left(\delta_\varphi(H) +1 \right)}{n-1} + \delta_\varphi(H)-\delta_\psi(u^2) \right]\hat{H}\\
	&~~~~ +\frac{\left(\delta_\varphi(H) +1 \right)^2}{n-1}f_{11}^2 + \frac{2\left(\delta_\varphi(H) +1 \right)f_{11}\psi(u^2)}{(n-1)\varphi(H)} + \frac{\psi(u^2)^2}{\varphi(H)^2}.
	\end{aligned}
	\end{align*}
\end{proof}

Motivated by \cite{he2023gradient}, we consider a weighted linearization operator
\begin{align}\label{weighted_linear}
\tilde{\mathcal{L}}_\varphi(\eta) := \mathcal{W}(\eta)^{-1} \operatorname{div}\left(\mathcal{W}(\eta)\varphi(H)\mathcal{A}(\nabla \eta) \right).
\end{align}
By Lemma \ref{LH_lower}, direct calculation gives
\begin{align}\label{weighted_linear2}
\begin{aligned}
\tilde{\mathcal{L}}_\varphi(\hat{H}) &\geqslant 2\varphi(H)\operatorname{Ric}(\nabla f) + \varphi(H)\left(\delta(H)+1+\frac{\delta_\mathcal{W}(\hat{H})}{2}\right)\frac{|\nabla \hat{H}|^2}{\hat{H}} \\
&~~~~   + \left[\frac{2(\delta_{\varphi}(H)+1)^2}{n-1} - 4\delta_{\varphi}'(H)H\right]\varphi(H) \hat{H}^2\\
&~~~~ + \left[\frac{2(\delta_{\varphi}(H)+1)^2}{n-1}-2(\delta_{\varphi}(H)+1) -2\delta_{\varphi}'(H)H\right]\varphi(H)\left\langle \nabla \hat{H}, \nabla f\right\rangle\\
&~~~~ + \frac{2\varphi(H)}{n-1}\left((\delta_{\varphi}(H)+1)\frac{\left\langle\nabla \hat{H}, \nabla f\right\rangle}{2\hat{H}} + \frac{\psi(u^2)}{\varphi(H)} \right)^2 \\
&~~~~ + {\delta_\mathcal{W}(\hat{H})}\varphi'(H)H\left(\frac{\left\langle\nabla \hat{H}, \nabla f\right\rangle}{\hat{H}}\right)^2\\
&~~~~ + 2\psi(u^2)\left[\frac{2\left(\delta_\varphi(H) +1 \right)}{n-1} + \delta_\varphi(H)-\delta_\psi(u^2) \right]\hat{H}.
\end{aligned}
\end{align}

At the end of this section, we present the following Sobolev inequality on Riemannian manifolds, which is critical to run the iteration.
\begin{theorem}[\cite{SaloffCoste1992UniformlyEO}]\label{sblv_thm}
	For $n>2$,  let $(M^n,g)$ be a complete Riemannian $n$-manifold with Ricci curvature bounded from below by $\operatorname{Ric} \geqslant -K$ {for some} $K\geqslant0$, then there exists $C$, depending only on $n$, such that for ball $B(R) \subset M$ with radius $R$ and volume $V(R)$, we have for any $f \in C_0^{\infty}(B)$,
	$$
	\begin{gathered}\label{sblv_inq}
	\left(\int_B|f|^{2 \chi} \right)^{1 / \chi} \leq e^{C(1+\sqrt{K}R)} V^{-2 / n} R^2\left(\int_B\left(|\nabla f|^2+R^{-2}|f|^2\right) \right), \\
	\end{gathered}
	$$
	where $\chi=n /(n-2)$. Meanwhile, for $n \leq 2$, the above inequality holds with $n$ replaced by any fixed $n^{\prime}>2$.
\end{theorem}

\section{Proof of the main theorem}
Let $(M,g)$ be a complete Riemannian manifold and $u$ be a positive local solution over an open neighborhood $\Omega$ containing $o\in M_u$  (otherwise, if $o\in M \setminus M_u$, Theorem \ref{mainthm} holds naturally). 

Note the curvature condition $\operatorname{Ric} \geqslant -K$  and adopt the conditions ($\varphi_1$) and ($\varphi_2$) in Lemma \ref{LH_lower}. Thus, after ignoring some nonnegative terms on the RHS of (\ref{eq-LH_lower}), one may deduce that
\begin{align}\label{btf}
\begin{aligned}
\mathcal{L}_\varphi(\hat{H}) &\geqslant -2K\varphi(H)\hat{H} + \left[\frac{2(\delta(H)+1)^2}{n-1} - 4\delta'(H)H\right]\varphi(H) \hat{H}^2\\
&~~~~ + \left[\frac{2(\delta(H)+1)^2}{n-1}-2(\delta(H)+1) -2\delta'(H)H\right]\varphi(H)\left\langle \nabla \hat{H}, \nabla f\right\rangle\\
&~~~~ +
2\psi(u^2)\left[\frac{2\left(\delta_\varphi(H) +1 \right)}{n-1} + \delta_\varphi(H)-\delta_\psi(u^2) \right]\hat{H}\\
&\geqslant -2K\varphi(H)\hat{H} + 2\gamma_{\varphi}\varphi(H) \hat{H}^2 - a_0\varphi(H)|\nabla \hat{H}||\nabla f|\\
&~~~~ +
2\psi(u^2)\left[\frac{2\left(\delta_\varphi(H) +1 \right)}{n-1} + \delta_\varphi(H)-\delta_\psi(u^2) \right]\hat{H},
\end{aligned}
\end{align}
where
$$
a_0 := \varGamma_{\varphi} + (d_\varphi + 1)^2 + 2(d_\varphi + 1).
$$
Since (\ref{btf}) holds only on $M_\varepsilon$, it follows that
\begin{align}\label{wbtf}
\begin{aligned}
&\int_{\Omega\cap M_\varepsilon} \left\langle \varphi(H)\nabla \hat{H} + 2\varphi'(H) \left\langle \nabla u, \nabla \hat{H}   \right\rangle \nabla u, \nabla \phi \right\rangle
\\&\leqslant 2K\int_{\Omega\cap M_\varepsilon}\varphi(H)\hat{H}\phi - 2\gamma_{\varphi}\int_{\Omega\cap M_\varepsilon} \varphi(H) \hat{H}^2\phi +  a_0\int_{\Omega\cap M_\varepsilon}\varphi(H)|\nabla \hat{H}||\nabla f|\phi\\
&~~~~ -
2\int_{\Omega\cap M_\varepsilon}\psi(u^2)\left[\frac{2\left(\delta_\varphi(H) +1 \right)}{n-1} + \delta_\varphi(H)-\delta_\psi(u^2) \right]\hat{H}\phi,
\end{aligned}
\end{align}
for any nonnegative test function $\phi$ compactly supported in $\Omega\cap M_\varepsilon$.

For the same $\varepsilon >0$, we take $\hat{H}_\varepsilon := \left(\hat{H} -\varepsilon\right)^+$, so that $\hat{H}_\varepsilon$ is compactly supported in $M_\varepsilon$. Since $\varphi(H) > 0$ holds in $M_\varepsilon$, it is valid to choose test function
$$\phi :=\frac{\lambda\hat{H}_\varepsilon^{b}\eta^2}{\varphi(H)}, $$
where $\lambda(x)$ is the characteristic function of $\left\{x\in \Omega: u^2(x)\in I_\psi\right\}$, and the cutoff function $\eta\in C^\infty_0(\Omega)$ and constant $b>1$ will be determined later. 
Then the last term of (\ref{wbtf}) is non-positive due to 
$$\psi(u^2)\left[\frac{2\left(\delta_\varphi(H) +1 \right)}{n-1} + \delta_\varphi(H)-\delta_\psi(u^2) \right]\geqslant0$$
on $\{u^2\in I_\psi\}$.

The first derivative of this test function $\phi$ is
\begin{align*}
\begin{aligned}
\nabla\phi &= \frac{b\lambda\hat{H}_\varepsilon^{b-1}\eta^2}{\varphi(H)} \nabla\hat{H} + \frac{2\lambda\hat{H}_\varepsilon^b\eta}{\varphi(H)}\nabla \eta - \frac{\lambda\varphi'(H)\hat{H}_\varepsilon^{b}\eta^2}{\varphi(H)^2}\nabla H \\
&= \left(\frac{b}{\varphi(H)}\hat{H}_\varepsilon^{b-1}\lambda - \frac{\varphi'(H)H}{\varphi(H)^2}\frac{\hat{H}_\varepsilon^{b}\lambda}{\hat{H}} \right)\eta^2\nabla\hat{H} + \frac{2\hat{H}_\varepsilon^b\eta\lambda}{\varphi(H)}\nabla \eta -\frac{2\varphi'(H)H}{\varphi(H)^2}{\hat{H}_\varepsilon^{b}\lambda}\nabla f.
\end{aligned}
\end{align*}
Thus, the LHS of (\ref{wbtf}) is then equal to
\begin{align}\label{wbtf_lhs}
\begin{aligned}
\int_{\Omega\cap M_\varepsilon}  &\left\langle \varphi(H)\nabla \hat{H} + 2\varphi'(H) \left\langle \nabla u, \nabla \hat{H}   \right\rangle \nabla u, \nabla \phi \right\rangle
\\
=&\int_{\Omega\cap M_\varepsilon}\left({b}\hat{H}_\varepsilon^{b-1} - \frac{\varphi'(H)H}{\varphi(H)}\frac{\hat{H}_\varepsilon^{b}}{\hat{H}} \right)\left(|\nabla \hat{H}|^2 + \frac{2\varphi'(H)}{\varphi(H)} \left\langle \nabla u, \nabla \hat{H}\right\rangle^2 \right)\lambda\eta^2\\
&+2\int_{\Omega\cap M_\varepsilon}\hat{H}_\varepsilon^b \left(\left\langle \nabla \hat{H}, \nabla \eta\right\rangle + \frac{2\varphi'(H)}{\varphi(H)}\left\langle \nabla u, \nabla \hat{H}\right\rangle\left\langle \nabla u, \nabla \eta\right\rangle\right)\lambda\eta\\
&-2\int_{\Omega\cap M_\varepsilon}\frac{\varphi'(H)H}{\varphi(H)}\hat{H}_\varepsilon^{b}\left(\left\langle \nabla f, \nabla \hat{H}\right\rangle + \frac{2\varphi'(H)}{\varphi(H)} \left\langle \nabla u, \nabla \hat{H}\right\rangle\left\langle \nabla u, \nabla f\right\rangle\right)\lambda\eta^2.
\end{aligned}
\end{align}
Noting that $\hat{H}_\epsilon \leqslant \hat{H}$ and
$$
\frac{\varphi'(H)H}{\varphi(H)} \leqslant \frac{d_\varphi}{2},
$$
we can observe the first term on the RHS of (\ref{wbtf_lhs}) could be estimated from below by
\begin{align}\label{wbtf_lhs1}
\begin{aligned}
\int_{\Omega\cap M_\varepsilon}&\left({b}\hat{H}_\varepsilon^{b-1} - \frac{\varphi'(H)H}{\varphi(H)}\frac{\hat{H}_\varepsilon^{b}}{\hat{H}} \right)\left(|\nabla \hat{H}|^2 + \frac{2\varphi'(H)}{\varphi(H)} \left\langle \nabla u, \nabla \hat{H}\right\rangle^2 \right)\lambda\eta^2\\
\geqslant& \int_{\{\varphi'(H) < 0\}}b\hat{H}_\varepsilon^{b-1}\left(|\nabla \hat{H}|^2 + l_\varphi|\nabla \hat{H}|^2 \right)\lambda\eta^2 + \int_{\{\varphi'(H) \geqslant 0\}}\left({b} - \frac{d_\varphi}{2}\right)\hat{H}_\varepsilon^{b-1}|\nabla \hat{H}|^2 \lambda\eta^2\\
\geqslant& \frac{a_1 b}{2} \int_{\Omega\cap M_\varepsilon}\hat{H}_\varepsilon^{b-1}|\nabla \hat{H}|^2\lambda\eta^2,
\end{aligned}
\end{align}
for $b$ is large enough (\textit{i.e.} $b>d_\varphi$) and
$$a_1 := \min\{1+l_\varphi,1\}.$$
Moreover, the second term on the RHS of (\ref{wbtf_lhs}) becomes 
\begin{align}\label{wbtf_lhs2}
\begin{aligned}
&2\int_{\Omega\cap M_\varepsilon}\hat{H}_\varepsilon^b \left(\left\langle \nabla \hat{H}, \nabla \eta\right\rangle + \frac{2\varphi'(H)}{\varphi(H)}\left\langle \nabla u, \nabla \hat{H}\right\rangle\left\langle \nabla u, \nabla \eta\right\rangle\right)\lambda\eta \\
\geqslant& -2\int_{\Omega\cap M_\varepsilon}\hat{H}_\varepsilon^b \left(|\nabla \hat{H}|| \nabla \eta| + \left|\frac{2\varphi'(H)}{\varphi(H)}\right||\nabla u|^2 |\nabla \hat{H}||\nabla \eta|\right)\lambda\eta\\
\geqslant&-a_2 \int_{\Omega\cap M_\varepsilon}\hat{H}_\varepsilon^b |\nabla \hat{H}|| \nabla \eta|\lambda\eta,
\end{aligned}
\end{align}
by setting
$$a_2 := 4\max \left\{1, |l_\varphi|, |d_\varphi| \right\}.$$
Finally the last one on the RHS of (\ref{wbtf_lhs}) turns into
\begin{align}\label{wbtf_lhs3}
\begin{aligned}
&-2\int_{\Omega\cap M_\varepsilon}\frac{\varphi'(H)H}{\varphi(H)}\hat{H}_\varepsilon^{b}\left(\left\langle \nabla f, \nabla \hat{H}\right\rangle + \frac{2\varphi'(H)}{\varphi(H)} \left\langle \nabla u, \nabla \hat{H}\right\rangle\left\langle \nabla u, \nabla f\right\rangle\right)\lambda\eta^2\\
\geqslant& -\int_{\Omega\cap M_\varepsilon}\left|\frac{\varphi'(H)H}{\varphi(H)}\right|\hat{H}_\varepsilon^{b}\left(|\nabla f||\nabla \hat{H}| + \left|\frac{2\varphi'(H)H}{\varphi(H)}\right| |\nabla f||\nabla \hat{H}|\right)\lambda\eta^2\\
\geqslant& -a_3\int_{\Omega\cap M_\varepsilon}\hat{H}_\varepsilon^{b}|\nabla f||\nabla \hat{H}|\lambda\eta^2,
\end{aligned}
\end{align}
where it could be chosen by
$$a_3 = {\max\left\{|l_\varphi|,|d_\varphi|\right\}}\cdot \max \left\{1, |l_\varphi|, |d_\varphi| \right\}. $$

After combining (\ref{wbtf_lhs1}), (\ref{wbtf_lhs2}) (\ref{wbtf_lhs3}) with (\ref{wbtf_lhs}), then (\ref{wbtf}) leads to
\begin{align*}
\begin{aligned}
\frac{a_1 b}{2} \int_{\Omega\cap M_\varepsilon}&\hat{H}_\varepsilon^{b-1}|\nabla \hat{H}|^2\lambda\eta^2 + 2\gamma_{\varphi}\int_{\Omega\cap M_\varepsilon} \hat{H}^2\hat{H}_\varepsilon^{b}\lambda\eta^2\\ \leqslant& 2K\int_{\Omega\cap M_\varepsilon}\hat{H}\hat{H}_\varepsilon^{b}\lambda\eta^2  +  (a_0 + a_3)\int_{\Omega\cap M_\varepsilon}\hat{H}_\varepsilon^{b}|\nabla \hat{H}||\nabla f|\lambda\eta^2
+ a_2\int_{\Omega\cap M_\varepsilon}\hat{H}_\varepsilon^{b}|\nabla \hat{H}||\nabla \eta|\lambda\eta
\\ \leqslant& 2K\int_{\Omega\cap M_\varepsilon}\hat{H}^{b+1}\lambda\eta^2  +  (a_0 + a_3)\int_{\Omega\cap M_\varepsilon}\hat{H}^{b}|\nabla \hat{H}||\nabla f|\lambda\eta^2
+ a_2\int_{\Omega\cap M_\varepsilon}\hat{H}^{b}|\nabla \hat{H}||\nabla \eta|\lambda\eta.
\end{aligned}.
\end{align*}
The last inequality is because $\hat{H}_\varepsilon \leqslant \hat{H}$. Note that $\hat{H}_\varepsilon\eta$ has compact support in $ M_\varepsilon\cap\Omega$, the integral can be extended to $\Omega$. Then by Fatou's lemma, we obtain that
\begin{align}\label{int_esti_1}
\begin{aligned}
\frac{a_1 b}{2} \int_{\Omega}&\hat{H}^{b-1}|\nabla \hat{H}|^2\lambda\eta^2 + 2\gamma_{\varphi}\int_{\Omega} \hat{H}^{b+2}\lambda\eta^2\\ \leqslant&\lim\limits_{\varepsilon \to 0}\frac{a_1 b}{2} \int_{\Omega}\hat{H}_\varepsilon^{b-1}|\nabla \hat{H}|^2\lambda\eta^2 + \lim\limits_{\varepsilon \to 0}2\gamma_{\varphi}\int_{\Omega} \hat{H}^2\hat{H}_\varepsilon^{b}\lambda\eta^2
\\ \leqslant& 2K\int_{\Omega}\hat{H}^{b+1}\lambda\eta^2  +  (a_0 + a_3)\int_{\Omega}\hat{H}^{b}|\nabla \hat{H}||\nabla f|\lambda\eta^2
+ a_2\int_{\Omega}\hat{H}^{b}|\nabla \hat{H}||\nabla \eta|\lambda\eta.
\end{aligned}
\end{align}

Again by mean of the {Cauchy's} inequality, the last two terms on the RHS of (\ref{int_esti_1}) could be estimated, respectively, as follows
\begin{align*}
\begin{aligned}
(a_0 + a_3)\int_{\Omega}\hat{H}^{b}|\nabla \hat{H}||\nabla f|\lambda\eta^2 \leqslant \frac{(a_0 + a_3)^2}{4\gamma_{\varphi}}\int_{\Omega} \hat{H}^{b-1}|\nabla \hat{H}|^2\lambda\eta^2+ \gamma_{\varphi}\int_{\Omega} \hat{H}^{b+2}\lambda\eta^2,
\end{aligned}
\end{align*}
and
\begin{align*}
\begin{aligned}
a_2\int_{\Omega}\hat{H}^{b}|\nabla \hat{H}||\nabla \eta|\lambda\eta \leqslant \frac{a_1b}{4}\int_{\Omega} \hat{H}^{b-1}|\nabla \hat{H}|^2\lambda\eta^2+ \frac{a_2^2}{a_1b}\int_{\Omega} \hat{H}^{b+1}|\nabla \eta|^2\lambda.
\end{aligned}
\end{align*}
We additionally requiring
\begin{align}\label{require_b}
b>\max\left\{\frac{2(a_0 + a_3)^2}{a_1\gamma_{\varphi}},d_{\varphi},1\right\},
\end{align}
so that (\ref{int_esti_1}) becomes
\begin{align}\label{int_esti_2}
\begin{aligned}
\frac{a_1 b}{8} \int_{\Omega}  \hat{H}^{b-1}|\nabla \hat{H}|^2\lambda\eta^2 + \gamma_{\varphi}\int_{\Omega} \hat{H}^{b+2}\lambda\eta^2
\leqslant 2K\int_{\Omega} \hat{H}^{b+1}\lambda\eta^2
+ \frac{a_2^2}{a_1b}\int_{\Omega} \hat{H}^{b+1}|\nabla \eta|^2\lambda.
\end{aligned}
\end{align}

Since $\lambda^s = \lambda$ for any $s>0$, from the inequality that
\begin{align*}
\left|\nabla\left(\hat{H}^{b/2+1/2}\eta\right)\right|^2 &\leqslant \frac{1}{2}\left(b+1\right)^2\hat{H}^{b-1}|\nabla \hat{H}|^2\eta^2 + 2\hat{H}^{b+1}|\nabla \eta|^2 \\
&\leqslant 2b^2\hat{H}^{b-1}|\nabla \hat{H}|^2\eta^2 + 2\hat{H}^{b+1}|\nabla \eta|^2,
\end{align*}
we get
\begin{align}\label{int_esti}
\begin{aligned}
\int_{\Omega} \left|\nabla\left(\hat{H}^{b/2+1/2}\eta\right)\right|^2\lambda + &\frac{16\gamma_\varphi b}{a_1} \int_{\Omega} \hat{H}^{ b+2}\lambda\eta^2\\
\leqslant& \frac{32Kb}{a_1}\int_{\Omega} \hat{H}^{b+1}\lambda\eta^2
+ \frac{16(a_2^2 + a_1^2)}{8a_1^2}\int_{\Omega} \hat{H}^{b+1}|\nabla \eta|^2\lambda.
\end{aligned}
\end{align}

When it comes to case that $u^2(x) \notin I_{\psi}$, we need the weighted operator in (\ref{weighted_linear}) and set $\mathcal{W}(\eta) = \eta^\alpha$ for some $\alpha>0$ which will be determined later, Then $\delta_\mathcal{W} \equiv 2\alpha$ and (\ref{weighted_linear2}) becomes
\begin{align*}
\begin{aligned}
\tilde{\mathcal{L}}_\varphi(\hat{H}) &\geqslant -2K\varphi(H)\hat{H}   + 2\gamma_{\varphi}\varphi(H)\hat{H}^2 + \varphi(H)\left(\delta_{\varphi}(H)+1+{\alpha}\right)\frac{|\nabla \hat{H}|^2}{\hat{H}}\\
&~~~~ - a_0\varphi(H)|\nabla \hat{H}||\nabla f|+\varphi(H)\left[\frac{(\delta_{\varphi}(H)+1)^2}{2(n-1)} + {\alpha}\delta_{\varphi}(H)\right]\frac{\left\langle\nabla \hat{H}, \nabla f\right\rangle^2}{\hat{H}^2}\\
&~~~~  + \frac{2(\delta_{\varphi}(H)+1)}{(n-1)}\cdot\frac{\left\langle\nabla \hat{H}, \nabla f\right\rangle}{\hat{H}}\psi(u^2)+ \frac{2}{n-1}\cdot\frac{\psi(u^2)^2}{\varphi(H)}\\
&~~~~  + 2\left[\frac{2\left(\delta_\varphi(H) +1 \right)}{n-1} + \delta_\varphi(H)-\delta_\psi(u^2) \right]\hat{H}\psi(u^2).
\end{aligned}
\end{align*}
According to $\delta_{\varphi}+1>0$ and
$$\left\langle\nabla \hat{H}, \nabla f\right\rangle^2\leqslant|\nabla\hat{H}|^2|\nabla f|^2 = |\nabla\hat{H}|^2\hat{H},$$
it follows that
\begin{align*}
\begin{aligned}
\tilde{\mathcal{L}}_\varphi(\hat{H}) &\geqslant -2K\varphi(H)\hat{H}   + 2\gamma_{\varphi}\varphi(H)\hat{H}^2 - a_0\varphi(H)|\nabla \hat{H}||\nabla f|\\
&~~~~ +\varphi(H)\left[\frac{(\delta_{\varphi}(H)+1)^2}{2(n-1)} + \alpha(\delta_{\varphi}(H)+1)\right]\frac{\left\langle\nabla \hat{H}, \nabla f\right\rangle^2}{\hat{H}^2}\\
&~~~~  + \frac{2(\delta_{\varphi}(H)+1)}{(n-1)}\cdot\frac{\left\langle\nabla \hat{H}, \nabla f\right\rangle}{\hat{H}}\psi(u^2)+ \frac{2}{n-1}\cdot\frac{\psi(u^2)^2}{\varphi(H)}\\
&~~~~  + 2\left[\frac{2\left(\delta_\varphi(H) +1 \right)}{n-1} + \delta_\varphi(H)-\delta_\psi(u^2) \right]\hat{H}\psi(u^2).
\end{aligned}
\end{align*}
Then by using $x^2 + 2xy \geqslant -y^2$ twice, we have
\begin{align*}
\begin{aligned}
\frac{\tilde{\mathcal{L}}_\varphi(\hat{H})}{\varphi(H)}
&\geqslant -2K\hat{H}   + 2\gamma_{\varphi}\hat{H}^2- a_0|\nabla \hat{H}||\nabla f|\\
&~~~~ +\left[\frac{2}{n-1} - \frac{2(\delta_{\varphi}(H)+1)^2}{2\alpha(\delta_{\varphi}(H)+1)(n-1)^2+(\delta_{\varphi}(H)+1)^2(n-1) }\right]\frac{\psi(u^2)^2}{\varphi(H)^2} \\
&~~~~ +  2\left[\frac{2\left(\delta_\varphi(H) +1 \right)}{n-1} + \delta_\varphi(H)-\delta_\psi(u^2) \right]\hat{H}\cdot\frac{\psi(u^2)}{\varphi(H)}\\
&\geqslant -2K\hat{H} - a_0|\nabla \hat{H}||\nabla f|\\
&~~~~ +2\left[\gamma_{\varphi}-\left(\frac{2\left(\delta_\varphi(H) +1 \right)}{n-1} + \delta_\varphi(H)-\delta_\psi(u^2)\right)^2\left(\frac{n-1}{4}+\frac{(\delta_{\varphi}(H)+1)}{8\alpha}\right)\right]\hat{H}^2\\
&\geqslant -2K\hat{H} - a_0|\nabla \hat{H}||\nabla f|\\
&~~~~ +2\left[\gamma_{\varphi}-\left(\frac{2\left(\delta_\varphi(H) +1 \right)}{n-1} + \delta_\varphi(H)-\delta_\psi(u^2)\right)^2\left(\frac{n-1}{4}+\frac{(d_{\varphi}+1)}{8\alpha}\right)\right]\hat{H}^2.
\end{aligned}
\end{align*}
{So in the weak sense, }
\begin{align}\label{wbtf2}
\begin{aligned}
&\int_{M_u} \left\langle \varphi(H)\hat{H}^\alpha\nabla \hat{H} + 2\varphi'(H)\hat{H}^\alpha \left\langle \nabla u, \nabla \hat{H}   \right\rangle \nabla u, \nabla \phi \right\rangle
\\\leqslant& 2K\int_{M_u}\varphi(H)\hat{H}^{1+\alpha}\phi  +  a_0\int_{M_u}\varphi(H)|\nabla \hat{H}||\nabla f|\hat{H}^{\alpha}\phi\\
&-
2\int_{M_u}\left[\gamma_{\varphi}-\left(\frac{2\left(\delta_\varphi(H) +1 \right)}{n-1} + \delta_\varphi(H)-\delta_\psi(u^2)\right)^2\left(\frac{n-1}{4}+\frac{(d_{\varphi}+1)}{8\alpha}\right)\right]\varphi(H)\hat{H}^{2+\alpha}\phi,
\end{aligned}
\end{align}

Similarly, choose the test function as
$$\phi :=\frac{\bar{\lambda}\hat{H}_\varepsilon^{b-\alpha}\eta^2}{\varphi(H)} $$
where $\bar{\lambda}$ is the characteristic function of $\left\{x\in \Omega: u^2(x)\notin I_\psi\right\}$, then the condition (\ref{psi2}) infers that the constant
$$\theta(\gamma_{\varphi},\varTheta_{\varphi,\psi}):= \gamma_{\varphi}-\sup_{\substack{s>0,\\t\in \mathbb{R}^+ -  I_{\psi}}} \left(\frac{2\left(\delta_\varphi(s) +1 \right)}{n-1} + \delta_\varphi(s)-\delta_\psi(t)\right)^2\frac{n-1}{4} >0.$$
Therefore, there exists a positive constant 
$$\alpha(\gamma_{\varphi},d_{\varphi},\varTheta_{\varphi,\psi}) :=\frac{1}{4\theta}\sup_{\substack{s>0,\\t\in \mathbb{R}^+ -  I_{\psi}}} \left(\frac{2\left(\delta_\varphi(s) +1 \right)}{n-1} + \delta_\varphi(s)-\delta_\psi(t)\right)^2(d_\varphi+1)>0,$$ 
such that the last term on the RHS of (\ref{wbtf2}) could be estimated from below as
\begin{align*}
\begin{aligned}
&2\int_{M_u}\left[\gamma_{\varphi}-\left(\frac{2\left(\delta_\varphi(H) +1 \right)}{n-1} + \delta_\varphi(H)-\delta_\psi(u^2)\right)^2\left(\frac{n-1}{4}+\frac{(d_{\varphi}+1)}{8\alpha}\right)\right]\varphi(H)\hat{H}^{2+\alpha}\phi\\
=&2\int_{\{u^2(x)\notin I_{\psi}\}}\left[\gamma_{\varphi}-\left(\frac{2\left(\delta_\varphi(H) +1 \right)}{n-1} + \delta_\varphi(H)-\delta_\psi(u^2)\right)^2\left(\frac{n-1}{4}+\frac{(d_{\varphi}+1)}{8\alpha}\right)\right]\varphi(H)\hat{H}^{2+\alpha}\phi\\
\geqslant&2\int_{\{u^2(x)\notin I_{\psi}\}}\left[\gamma_{\varphi}-\sup_{\substack{s>0,\\t\in \mathbb{R}^+ -  I_{\psi}}}\left(\frac{2\left(\delta_\varphi(s) +1 \right)}{n-1} + \delta_\varphi(s)-\delta_\psi(t)\right)^2\left(\frac{n-1}{4}+\frac{(d_{\varphi}+1)}{8\alpha}\right)\right]\varphi(H)\hat{H}^{2+\alpha}\phi\\
\geqslant&2\int_{\{u^2(x)\notin I_{\psi}\}}\left[\theta-\sup_{\substack{s>0,\\t\in \mathbb{R}^+ -  I_{\psi}}}\left(\frac{2\left(\delta_\varphi(s) +1 \right)}{n-1} + \delta_\varphi(s)-\delta_\psi(t)\right)^2\left(\frac{d_{\varphi}+1}{8\alpha}\right)\right]\varphi(H)\hat{H}^{2+\alpha}\phi\\
=&\theta\int_{M_u}\hat{H}^{2+\alpha}\hat{H}_\varepsilon^{b-\alpha}\eta^2\bar{\lambda}.
\end{aligned}
\end{align*}
Following the same process from (\ref{wbtf_lhs}) to (\ref{int_esti}), we obtain
\begin{align}\label{int_esti_c}
\begin{aligned}
\int_{\Omega} \left|\nabla\left(\hat{H}^{b/2+1/2}\eta\right)\right|^2\bar{\lambda} + &a_4\theta b \int_{\Omega} \hat{H}^{ b+2}\bar{\lambda}\eta^2\\
\leqslant& a_5Kb\int_{\Omega} \hat{H}^{b+1}\bar{\lambda}\eta^2
+ a_6\int_{\Omega} \hat{H}^{b+1}|\nabla \eta|^2\bar{\lambda},
\end{aligned}
\end{align}
for constant 
\begin{align}\label{require_b2}
b>\max\left\{\frac{a_7}{\theta},d_{\varphi},\alpha\right\}.
\end{align}
Noticing that $\lambda + \bar{\lambda} \equiv 1$, one may deduce from (\ref{int_esti}) and (\ref{int_esti_c}) that
\begin{align}\label{int_esti_f}
\begin{aligned}
\int_{\Omega} \left|\nabla\left(\hat{H}^{b/2+1/2}\eta\right)\right|^2 + &a_7 b \int_{\Omega} \hat{H}^{ b+2}\eta^2\\
\leqslant& a_8Kb\int_{\Omega} \hat{H}^{b+1}\eta^2
+ a_9\int_{\Omega} \hat{H}^{b+1}|\nabla \eta|^2,
\end{aligned}
\end{align}
by adjusting the coefficients if necessary.

Then if let $\Omega = B(o,R)$, Theorem \ref{sblv_thm} shows that when $n>2$
\begin{align}
\begin{aligned}\label{int_esti_sblv}
\left(\int_{\Omega} \hat{H}^{(b+1)\chi} \eta^{2\chi}\right)^{1/\chi}\leqslant e^{C(1+\sqrt{K}R)}V^{-2/n}\left(R^2\int_{\Omega}\left|\nabla\left(\hat{H}^{b/2+1/2}\eta\right)\right|^2+ \int_{\Omega} \hat{H}^{b+1} \eta^{2}\right).
\end{aligned}
\end{align}
with $\chi=n/(n-2)$. Set $b_0 = c_0\left(1+\sqrt{K}R\right)$ and choose $c_0$ large enough to satisfy (\ref{require_b}) and (\ref{require_b2}), in which $b = b_0$ may be determined later. In terms of (\ref{int_esti}) and (\ref{int_esti_sblv}), then direct calculation implies that
\begin{align*}
\begin{aligned}
&\left(\int_{\Omega} \hat{H}^{(b+1)\chi} \eta^{2\chi}\right)^{1/\chi} + a_7be^{c_1b_0}\left(\frac{R^2}{V^{2/n}}\right)\int_{\Omega} \hat{H}^{b+2}\eta^2\\
&~~~~~~~~~~~~~\leqslant (a_{8}bKR^2b+1)e^{c_1b_0}V^{-2/n}\int_{\Omega} \hat{H}^{b+1}\eta^{2} + a_{9}e^{c_1b_0}\left(\frac{R^2}{V^{2/n}}\right)\int_{\Omega} \hat{H}^{b+1}|\nabla\eta|^2.
\end{aligned}
\end{align*}
Noticing that
$$a_{8}bKR^2b+1 \leqslant \max\left\{{a_8},1\right\}\cdot(KR^2 + 1)b \leqslant \max\left\{{a_8},1\right\}\cdot(\sqrt{K}R + 1)^2b,$$
we have
\begin{align}\label{int_esti_final}
\begin{aligned}
&\left(\int_{\Omega} \hat{H}^{(b+1)\chi} \eta^{2\chi}\right)^{1/\chi} + a_7be^{c_1b_0}\left(\frac{R^2}{V^{2/n}}\right)\int_{\Omega} \hat{H}^{b+2}\eta^2\\
&~~~~~~~~~~~~~\leqslant a_{10}b_0^2be^{c_1b_0}V^{-2/n}\int_{\Omega} \hat{H}^{b+1}\eta^{2} + a_{10}e^{c_1b_0}\left(\frac{R^2}{V^{2/n}}\right)\int_{\Omega} \hat{H}^{b+1}|\nabla\eta|^2,
\end{aligned}
\end{align}
in which the constants could be chosen as
$$
c_1 = \frac{C}{c_0},~~~~a_{10}=\max\left\{\frac{\max\left\{{a_8},1\right\}}{c_0^2},a_{9}\right\}.
$$

Next we shall show an $L^{\beta}$ estimate for $\hat{H}$ as an initiate value of the iteration.
\begin{lemma}\label{Lb_esti}
	Under the same conditions above, take
	\begin{align}
	\beta_0 = b_0+1\text{ and } \beta_1 = \beta_0\chi.
	\end{align}
	Then there exists $a_{11}>0$ such that
	\begin{align}
	||\hat{H}||_{L^{\beta_1}(B_{3R/4})} \leqslant a_{11}\frac{b^2_0}{R^2}V^{1/\beta_1}.
	\end{align}
\end{lemma}
\begin{proof}
	Set $b = b_0$ and decompose the first term on RHS of (\ref{int_esti_final}) into two parts as
	\begin{align*}
	\begin{aligned}
	a_{10}b_0^3e^{c_1b_0}V^{-2/n}\int_{B_R} \hat{H}^{\beta_0}\eta^{2} = a_{10}b_0^3e^{c_1b_0}V^{-2/n}\left(\int_{\Omega_1} \hat{H}^{\beta_0}\eta^{2} + \int_{\Omega_2} \hat{H}^{\beta_0}\eta^{2}\right),
	\end{aligned}
	\end{align*}
	where
	\begin{align*}
	\Omega_1:= \left\{\hat{H}>\frac{2a_{10}b_0^2}{a_7R^2}\right\} \text{ and } \Omega_2:= \left\{\hat{H}\leqslant\frac{2a_{10}b_0^2}{a_7R^2}\right\}.
	\end{align*}
	This yields
	\begin{align}\label{RHS1_esti}
	\begin{aligned}
	a_{10}b_0^3e^{c_1b_0}V^{-2/n}&\int_{B_R} \hat{H}^{\beta_0}\eta^{2} \\
	\leqslant& \frac{1}{2}a_7b_0e^{c_1b_0}\left(\frac{R^2}{V^{2/n}}\right)\int_{B_R} \hat{H}^{\beta_0+1}\eta^2 +
	a_{10}a_{12}^{\beta_0}b_0^3e^{c_1b_0}V^{1-2/n}\left(\frac{b_0}{R}\right)^{2\beta_0},
	\end{aligned}
	\end{align}
	for $a_{12} = 2a_{10}/a_7$.
	
	To deal with the second term on the RHS of (\ref{int_esti_final}), we should choose the cutoff function $\eta = \eta_0^{\beta_0+1}$, where $\eta_0$ is a smooth function with compact support in $B(R)$ with $0\leqslant \eta_0 \leqslant 1$ and $\eta_0 \equiv 1$ in $B{(3R/4)}$, as well as satisfies that
	$$
	|\nabla \eta_0|\leqslant\frac{c_2(n)}{R}.
	$$
	Therefore,
	\begin{align}\label{eta_esti}
	|\nabla \eta|^2\leqslant c^2_2\left(\frac{\beta_0+1}{R} \right)^2\eta_0^{2\beta_0} = c^2_2\left(\frac{\beta_0+1}{R}\right)^2\eta^{\frac{2\beta_0}{\beta_0+1}}.
	\end{align}
	Substituting (\ref{eta_esti}) for the second term on the RHS of (\ref{int_esti_final}), we have
	\begin{align}\label{RHS2_esti}
	\begin{aligned}
	a_{10}b_0&  e^{c_1b_0}\left(\frac{R^2}{V^{2/n}}\right)\int_{B_R} \hat{H}^{\beta_0}|\nabla\eta|^2 \leqslant a_{10}c_2^2\left(\beta_0+1\right)^2e^{c_1b_0}V^{-2/n}\int_{B(R)} \hat{H}^{\beta_0}\eta^{\frac{2\beta_0}{\beta_0+1}}\\
	\leqslant& a_{10}c_2^2e^{c_1b_0}\left(\frac{\left(\beta_0+1\right)^2}{V^{2/n}}\right)\left(\int_{B(R)} \hat{H}^{\beta_0+1}\eta^{2}\right)^{\frac{\beta_0}{\beta_0+1}}\left(\int_{B(R)} 1\right)^{\frac{1}{\beta_0+1}}\\
	\leqslant&\frac{1}{2}a_7b_0e^{c_1b_0}\left(\frac{R^2}{V^{2/n}}\right)\int_{B(R)} \hat{H}^{\beta_0+1}\eta^{2}+ 2a_{10}c_2^2\left(\frac{4a_{10}c_2^2}{a_7}\right)^{\beta_0}\beta_0e^{c_1b_0}V^{1-2/n}\left(\frac{\beta_0}{R}\right)^{2\beta_0},
	\end{aligned}
	\end{align}
	where we have utilized the H\"older's inequality and Young's inequality at the last two inequalities, respectively.

	It follows from (\ref{int_esti_final}), (\ref{RHS1_esti}), and (\ref{RHS2_esti})  that
	\begin{align}
	\begin{aligned}
	\left(\int_{B(R)} \hat{H}^{\beta_0\chi} \eta^{2\chi}\right)^{1/\chi} \leqslant a_7a_{8}^{\beta_0}\beta_0^3e^{c_1b_0}V^{1/\chi}\left(\frac{\beta_0}{R}\right)^{2\beta_0},
	\end{aligned}
	\end{align}
	which is exact (\ref{Lb_esti}) after taking $(1/\beta_0)$-root on both sides.
\end{proof}

Now it is ready to finish our main theorem.

\begin{proof}[Proof of Theorem \ref{mainthm}]
	Here we go back to (\ref{int_esti_final}) and dismiss the second nonnegative term on the LHS. It follows that
	\begin{align}\label{int_esti_iteration}
	\begin{aligned}
	\left(\int_{B(R)} \hat{H}^{(b+1)\chi} \eta^{2\chi}\right)^{1/\chi}
	\leqslant a_{10}\left(\frac{e^{c_1b_0}}{V^{2/n}}\right)\int_{B(R)} \left(\beta_0^2(b+1)\eta^{2}+ R^2|\nabla\eta|^2 \right)\hat{H}^{b+1},
	\end{aligned}
	\end{align}
	where $\beta_0$ are given in Lemma \ref{Lb_esti}.
	
	We now choose the sequences of $\beta_k$ and $R_k$ by
	\begin{align*}
	\beta_1 = \beta_0\chi, ~~\beta_2 = \beta_0\chi^2, ~~\cdots ,~~ \beta_k = \beta_0\chi^k, \cdots,\\
	R_1 = \frac{3R}{4}, ~~R_2 = \frac{9R}{16}, ~~\cdots , ~~R_k = \frac{R}{2} + \frac{R}{4^k}, ~~\cdots,
	\end{align*}
	{so that}
	\begin{align*}
	\beta_k \to +\infty \text{ and } R_k \to \frac{R}{2}
	\end{align*}
	as $k \to \infty$. Moreover, one could choose a sequence of cutoff functions $\eta_k$ such that
	\begin{align}
	\left\{\begin{aligned}
	\eta_k \equiv 1~~~~~~~~~~~~~~~~~~ &\text{ ~~~in } B(R_{k+1}),\\
	0 \leqslant \eta_k \leqslant 1 \text{ and } |\nabla\eta_k|\leqslant\frac{c_3(n)4^k}{R} &\text{ ~~~in } B(R_{k}) - B(R_{k+1}),\\
	\eta_k \equiv 0~~~~~~~~~~~~~~~~~~ &\text{ ~~~in } B(R) - B(R_{k+1}).
	\end{aligned}\right.
	\end{align}
	By letting $b = b_k$ in (\ref{int_esti_iteration}) with $$b_k+1 = \beta_k,$$ and noting that
	$$
	b_k < \beta_{k} = b_k + 1 \leqslant 2b_k,
	$$
	we have
	\begin{align*}
	\begin{aligned}
	\left(\int_{B(R_{k+1})} \hat{H}^{\beta_{k+1}}\right)^{1/\chi}
	&\leqslant \left(\frac{a_{10}e^{c_1b_0}}{V^{2/n}}\right)\left(\beta_0^3\chi^k+ c_3 16^k \right)\int_{B(R_k)} \hat{H}^{\beta_k}\\
	&\leqslant \left(\frac{a_{10}\left(\beta_0^3+c_3\right)e^{c_1b_0}}{V^{2/n}}\right)\left( 16 \right)^k\int_{B(R_k)} \hat{H}^{\beta_k},
	\end{aligned}
	\end{align*}
	namely,
	\begin{align}\label{int_esti_iteration_2}
	\begin{aligned}
	||\hat{H}||_{L^{\beta_{k+1}}(B(R_{k+1}))} \leqslant \left(\frac{a_{10}\left(\beta_0^3+c_3\right)e^{c_1b_0}}{V^{2/n}}\right)^{1/\beta_k}\left( 16 \right)^{k/\beta_k}||\hat{H}||_{L^{\beta_{k}}(B(R_{k}))}.
	\end{aligned}
	\end{align}
	
	Then iterating (\ref{int_esti_iteration_2}) from $k = 1$ leads to 
	\begin{align}
	\begin{aligned}
	||\hat{H}||_{L^{\infty}(B(R/2))} &\leqslant 16^{\sum_{k=1}^{\infty} k/\beta_k}\left(\frac{a_{10}\left(\beta_0^3+c_3\right)e^{c_1b_0}}{V^{2/n}}\right)^{\sum_{k=1}^{\infty} 1/\beta_k} ||\hat{H}||_{L^{\beta_{1}}(B(3R/4))}\\
	&\leqslant e^{c_1}\left(16^{\frac{n^2}{4}}a_{10}^{\frac{n}{2}} \right)^{\frac{1}{\beta_1}}\left( \beta_0^3+c_3\right)^{\frac{n}{2\beta_1}}V^{-\frac{1}{\beta_1}}||\hat{H}||_{L^{\beta_{1}}(B(3R/4))},
	\end{aligned}
	\end{align}
	by noticing 
	\begin{align}
	\sum_{k=1}^{\infty}\frac{1}{\beta_k} = \frac{n}{2\beta_1} \text{ ~and~ } \sum_{k=1}^{\infty}\frac{k}{\beta_k} = \frac{n^2}{4\beta_1}.
	\end{align}
	According to the boundedness of $f(x) = C^{1/x}$ and $g(x) = (x+C)^{1/x}$ for $x>1$, we can find a constant $a_{13}$ independent of $b$, {such that}
	$$
	a_{13} \geqslant e^{c_1}\left(a_{10}^{\frac{n}{2}}16^{\frac{n^2}{4}} \right)^{\frac{1}{\beta_1}}\left( \beta_0^3+c_3\right)^{\frac{n}{2\beta_1}}.
	$$

	Finally, we conclude from Lemma \ref{Lb_esti} that
	\begin{align}\label{final_result}
	||\hat{H}||_{L^{\infty}(B_{R/2})}\leqslant
	a_{11}a_{13}\frac{b_0^2}{R^2} =
	C(n,d_{\varphi},l_{\varphi},\gamma_{\varphi},\varGamma_{\varphi},\varTheta_{\varphi,\psi})\frac{\left(1+\sqrt{K}R\right)^2}{R^2}
	\end{align}
	which finishes the proof for $n>2$.
	
	If $n=2$, Theorem \ref{sblv_thm} asserts that
	\begin{align*}
	\begin{aligned}
	\left(\int_{\Omega} \hat{H}^{\frac{(b+1)m}{m-2}} \eta^{\frac{2m}{m-2}}\right)^{\frac{m-2}{m}}\leqslant e^{C(1+\sqrt{K}R)}V^{-2/m}\left(R^2\int_{\Omega}\left|\nabla\left(\hat{H}^{b/2+1/2}\eta\right)\right|^2+ \int_{\Omega} \hat{H}^{b+1} \eta^{2}\right)
	\end{aligned}
	\end{align*}
	holds for each $m>2$. In particular, one can take $m=4$ and it follows that
	\begin{align}\label{int_esti_sblv_n2}
	\begin{aligned}
	\left(\int_{\Omega} \hat{H}^{2(b+1)} \eta^{4}\right)^{\frac{1}{2}}\leqslant e^{C(1+\sqrt{K}R)}V^{-\frac{1}{2}}\left(R^2\int_{\Omega}\left|\nabla\left(\hat{H}^{b/2+1/2}\eta\right)\right|^2+ \int_{\Omega} \hat{H}^{b+1} \eta^{2}\right).
	\end{aligned}
	\end{align}
	Then the same procedure for $n>2$ can be applied to this case, which yields
	\begin{align*}
	||\hat{H}||_{L^{2\beta_0}(B_{3R/4})} \leqslant a_{14}V^{\frac{1}{2\beta_0}}\cdot\frac{b_0^2}{R^2},
	\end{align*}
	and
	\begin{align*}
	\begin{aligned}
	||\hat{H}||_{L^{\infty}(B(R/2))} \leqslant  a_{15}V^{-\frac{1}{2\beta_0}}||\hat{H}||_{L^{2\beta_0}(B(3R/4))}.
	\end{aligned}
	\end{align*}
	Hence (\ref{final_result}) is also valid for $n=2$.
\end{proof}

Then we give the direct applications of the gradient estimates.
\begin{proof}[Proof of Theorem \ref{harnack}] Under the same conditions in Theorem 1.1,
	let  $x$, $y\in B(R)$ be any two points with minimal geodesic $l$ connecting them.
	Then using the gradient estimate and the fact that $\operatorname{length}(l)\leqslant 2R$, we have
	\begin{align}
	\begin{aligned}
	\log u(x)- \log u(y) \leqslant\int_{l}|\nabla\log u|&\leqslant \int_{l}C\frac{\sqrt{K}R+1}{R}\\
	&\leqslant 2C(\sqrt{K}R+1).
	\end{aligned}
	\end{align}
	Therefore,
	$$
	u(x)  \leqslant e^{C(1+\sqrt{\kappa} R)}u(y) .
	$$
\end{proof}

Another significance of the Cheng--Yau gradient estimate is to derive the Liouville theorems for some differential equations on complete but non-compact manifolds.
\begin{proof}[Proof of Theorem \ref{liou}]
	When $K=0$ and $0< u\leqslant A$ is a bounded positive solution of 
	$$\Delta_\varphi(u) + \psi(u^2)u = 0, $$ letting $R \to \infty$, we see $|\nabla u|=0$. Consequently, $u$ must be a constant and $\psi(u^2)u = 0$. if  $\psi(t)\neq0$ for any $t>0$, then there is no such positive solution for this equation. Otherwise, $u^2$ shall be a positive root of $\psi(t) = 0$.
\end{proof}

\section{Applications and some remarks}
In this section we will apply Theorem 1.1 to several specific examples for $\varphi(t)$ and $\psi(t)$.

\begin{example}\label{ex_pq}
	Assume $\psi(t)\equiv t^{p/2-1} + t^{q/2-1}$ we get the well-known $(p,q)$-Laplacian
	$$\Delta_{p,q} u := \operatorname{div}\left(\left(|\nabla u|^{p-2} + |\nabla u|^{q-2} \right)\nabla u  \right)=\Delta_{p}u + \Delta_{q}u.$$	
	
	Without loss of generality, suppose that $p<q$, then, by direct calculation, we have 
	$$\delta_\varphi(t) = \frac{(p-2)t^{p/2-1}+(q-2)t^{q/2-1} }{t^{p/2-1}+t^{q/2-1}} = (q-2)-\frac{(q-p)}{1+t^{(q-p)/2}}$$
	$$d_\varphi = (q-2) \text{ and }l_\varphi = (p-2).$$
	The condition ($\varphi_1$) follows that $q>p>1$. Also the $\varGamma_{\varphi}$ in condition ($\varphi_2$) can be determined by
	\begin{align*}
	&\frac{\left(\delta_{\varphi}(t)+1\right)^2 }{n-1} - 2t\delta_{\varphi}'(t) \\ &=\frac{(p-1)^2t^{p-2}+(q-1)^2t^{q-2}+\left(2(p-1)(q-1)-(n-1)(p-q)^2\right)t^{(p+q)/2-2} }{(n-1)\left(t^{p-2}+t^{q-2}+t^{(p+q)/2-2}\right)}\\
	&\leqslant \frac{(q-1)^2}{n-1}=:\varGamma_{\varphi}.
	\end{align*}
	Now we need to find the sufficient and necessary condition for the existence of $\gamma_{\varphi}>0$ in condition (${\varphi_2}$). In other word, there shall exist some $0<\gamma\leqslant{(q-1)^2}$, 
	such that, for any $X\geqslant0$,
	\begin{align}\label{abpq}
	\left((p-1)^2-\gamma\right)+\left(2(p-1)(q-1)-(n-1)(p-q)^2-2\gamma\right)X + \left((q-1)^2-\gamma \right)X^{2}\geqslant 0.
	\end{align}
	(\ref{abpq}) holds if and only if there exists some $0<\gamma\leqslant{(q-2)^2}$ such that
	\begin{align}\label{abpq1}
	2(p-1)(q-1)-(n-1)(p-q)^2-2\gamma>0,
	\end{align}
	or
	\begin{align}\label{abpq2}
	\left(2(p-1)(q-1)-(n-1)(p-q)^2-2\gamma\right)^2-4\left((p-1)^2-\gamma\right)\left((q-1)^2-\gamma \right)<0.
	\end{align}
	Combining (\ref{abpq1}) and (\ref{abpq2}) we have
	\begin{align}\label{pqn}
		\frac{n-1}{4}<\frac{(p-1)(q-1)}{(p-q)^2},
	\end{align}
	and the desired $$\gamma_{\varphi} = \frac{4(p-1)(q-1)-(n-1)(q-p)^2}{4n}.$$
\end{example}

	{For fixed $n$, we can draw (\ref{pqn}) in terms of the coordinate $(p,q)$ (see the Figure 1 (A)), which shows that the admissible area is indeed between two straight lines. From another perspective, when fixing $p$ (the Figure 1 (B)), the closer $q$ is to the $p$, the higher dimension Cheng--Yau estimate holds for. }

\begin{figure}
	\begin{subfigure}{0.4\textwidth}
		\begin{tikzpicture}[scale=0.8]
		\draw[-stealth] (-0.5,0)--(5,0) node[below]{$q$};
		\draw[-stealth] (0,-0.5)--(0,5) node[left]{$p$};
		\draw (0,0) node [below left] {$0$};
		\foreach \i in {1}{\draw (\i,0)--node [below] {$1$}(\i,0.05);}
		\draw (0,1) node [left] {$1$};
		\draw[domain=1:3,smooth,variable=\t] plot ({2*\t-1},\t);
		\draw[domain=1:3,smooth] plot(\x,{2*\x-1});
		\filldraw [fill=gray!20] plot [domain=1:3,smooth]  (\x,{2*\x-1}) -- (5,3)node[above left] {$\frac{4(p-1)(q-1)}{(p-q)^2}>n-1$};
		\draw [white] plot [domain=3:5,smooth]  (\x,{8-\x});
		\draw[domain=0:1,style=dashed] plot (\x,1);
		\draw[domain=0:1,style=dashed,variable=\t] plot (1,\t);
		\end{tikzpicture} 
		\caption{The admissible area for $(p,q)$ with fixed $n$}
		\label{fig:subim1}
	\end{subfigure}
	\begin{subfigure}{0.4\textwidth}
		\begin{tikzpicture}[scale=0.8]
		\draw[-stealth] (-0.5,0)--(5,0) node[below]{$|p-q|$};
		\draw[-stealth] (0,-0.5)--(0,5) node[left]{$n$};
		\draw (0,0) node [below left] {$0$};
		\foreach \i in {1}{\draw (\i,0)--node [below] {$1$}(\i,0.05);}
		\draw (0.1,1) node [left] {$1$};
		\draw[domain=0.43:5,smooth]
		plot (\x,{0.5*(\x+1)/(\x^2)+1});
		\filldraw [fill=gray!20](0,4.87)--plot [domain=0.43:5,color=white](\x,{0.5*(\x+1)/(\x^2)+1})--(5,1)--(0,1);
		\draw[domain=0:5,color=gray!20] plot (\x,1);
		\draw[domain=0:5,style=dashed] plot (\x,1);
		\draw[domain=0:0.43,color=gray!20] plot (\x,4.87);
		\draw[domain=1:1.13,color=gray!20] plot (5,\x);
		\draw[domain=0:5,smooth] plot (0,\x);
		\end{tikzpicture}
		\caption{The admissible area for $n$ and $|p-q|$ with fixed $p$}
		\label{fig:subim2}
	\end{subfigure}
\caption{}
\end{figure}

\
\begin{example}\label{ex_pr}
	Assume $\varphi(t) =  \sum_{i=1}^{r}a_it^{p_i/2-1}$, which means
	$$\tilde{\Delta}_{p_1,...p_r}u:=\left(\sum_{i=1}^ra_i\Delta_{p_i}\right)u = \operatorname{div}\left(\sum_{i=1}^{r}a_i|\nabla u|^{p_i-2}\nabla u  \right),$$
	where one could assume $a_i>0$ and {$p_1<...<p_r$} without loss of generality. Then
	$$\delta_\varphi(t) = \frac{\sum_{i=1}^{r}a_i(p_i-2)t^{p_i/2-1}}{\sum_{i=1}^{r}a_it^{p_i/2-1}},$$
	$$d_\varphi = (p_r-2) \text{ and }l_\varphi = (p_1-2), $$
	and {the} condition ($\varphi_1$) yields that $p_1>1$.
	When it comes to {the}  condition ($\varphi_2$), we need to compute
	\begin{align}
	\begin{aligned}
	2t\delta_{\varphi}'(t) &=\frac{\left(\sum a_i(p_i-2)^2t^{p_i/2-1}\right)\left(\sum a_it^{p_i/2-1}\right) - \left(\sum a_i(p_i-2)t^{p_i/2-1}\right)^2}{\left(\sum a_it^{p_i/2-1}\right)^2}\\
	&=\frac{\sum_{i\neq j}\left( (p_i-2)^2-(p_i-2)(p_j-2)\right)a_ia_jt^{(p_i+p_j)/2-2}}{\sum a_ia_jt^{(p_i+p_j)/2-2}}.
	\end{aligned}
	\end{align}
	After dividing the summation into $i>j$ and $i<j$, then switching the index $i$ with $j$, it becomes 
	\begin{align*}
	2t\delta_{\varphi}'(t) &=\frac{\sum_{j>i}\left( (p_i-p_j)(p_i-2) + (p_j-p_i)(p_j-2)\right)a_ia_jt^{(p_i+p_j)/2-2}}{\sum a_ia_jt^{(p_i+p_j)/2-2}}\\&=\frac{\sum_{j>i}\left( (p_i-p_j)^2\right)a_ia_jt^{(p_i+p_j)/2-2}}{\sum a_ia_jt^{(p_i+p_j)/2-2}}
	\\&=\frac{\sum\left( (p_i-p_j)^2\right)a_ia_jt^{(p_i+p_j)/2-2}}{2\sum a_ia_jt^{(p_i+p_j)/2-2}}\leqslant\frac{(p_r-p_1)^2}{2}.
	\end{align*}
	Note that $$\frac{\left(p_1-1\right)^2 }{n-1}\leqslant \frac{\left(\delta_{\varphi}(t)+1\right)^2 }{n-1}\leqslant\frac{\left(p_r-1\right)^2 }{n-1},$$
	one can set $$\gamma_{\varphi} = \frac{(p_1-1)^2}{n-1}-\frac{(p_r-p_1)^2}{2} \text{ and } \varGamma_{\varphi} = \frac{(p_r-1)^2}{n-1},$$
	provided 
	{
	\begin{align}\label{nNc}
		\frac{(p_1-1)^2}{(p_r-p_1)^2}>\frac{n-1}{2}.
	\end{align}
}
\end{example}

Hence, {it is an interesting phenomenon that the upper and lower bounds of the degree function of weighted $(p_1,...,p_r)$-Laplacian $\tilde{\Delta}_{p_1,...,p_r}$ is independent of the weight $a_i$, and precisely,} it turns out that {they only depend} on the {maximum and minimum of $p_i$, if} $p_i$ are large enough or very close to each other, 
	 {namely, this property can be reflected in the constant 
	\begin{align}
		\mathcal{N}_{1}:= 2\left(\frac{\min\{p_i\} - 1}{\max\{p_i\}- \min\{p_i\}}\right)^2 +1,
	\end{align}
	called the \textbf{first critical dimension} of $\tilde{\Delta}_{p_1,...,p_r}$. Then (\ref{nNc}) implies that $n<\mathcal{N}_1$.}
	
	{Further, we define the \textbf{second critical dimension} by 
	  \begin{align}
	  \mathcal{N}_{2}:= \sqrt{2\mathcal{N}_1 +3} -2,
	  \end{align} 
  and it is easy to check that $\mathcal{N}_2 < \mathcal{N}_1$}.
	  
{In the rest of this section, we will show that these constants play an important role in determining the gradient estimate and the Liouville property of the weighted $(p_1,...,p_r)$-Laplacian equation when $n$ is bounded by different critical dimensions. In particular, if $\tilde{\Delta}_{p_1,...,p_r}$ reduces to $p$-Laplacian, so that $p_1=p_r$, then the critical dimensions are defined to be $\infty$, which means $n<\mathcal{N}_1$ and $n<\mathcal{N}_2$ for any dimension, so the dimension has little effect on the Liouville property.}

	 {Similarly, we can also draw the following Figure 2 with respect to $p_1$, $p_r$ and $n$.}

\begin{figure}[H]
	\begin{subfigure}{0.4\textwidth}
		\begin{tikzpicture}[scale=0.8]
		\draw[-stealth] (-0.5,0)--(5,0) node[below]{$p_1$};
		\draw[-stealth] (0,-0.5)--(0,5) node[left]{$p_r$};
		\draw (0,0) node [below left] {$0$};
		\foreach \i in {1}{\draw (\i,0)--node [below] {$1$}(\i,0.05);}
		\draw (0,1) node [left] {$1$};
		\draw (4.5,2) node [above] {$\frac{(p_1-1)^2}{(p_r-q_1)^2}>n-1$};
		\draw[domain=0:5,style=dashed,variable=\t] plot (\t,\t);
		\draw[domain=1:3,smooth] plot(\x,{2*\x-1});
		\filldraw [fill=gray!20] plot [domain=1:3,smooth]  (\x,{2*\x-1}) -- (5,5);
		\draw [domain=3:5,smooth,color=white] plot   (\x,5);
		\draw[domain=0:1,style=dashed] plot (\x,1);
		\draw[domain=0:1,style=dashed,variable=\t] plot (1,\t);
		\end{tikzpicture} 
		\caption{The admissible area for $(p_1,p_r)$}
	\end{subfigure}
	\begin{subfigure}{0.4\textwidth}
		\begin{tikzpicture}[scale=0.8]
		\draw[-stealth] (-0.5,0)--(5,0) node[below]{$p_r-p_1$};
		\draw[-stealth] (0,-0.5)--(0,5) node[left]{$n$};
		\draw (0,0) node [below left] {$0$};
		\foreach \i in {1}{\draw (\i,0)--node [below] {$1$}(\i,0.05);}
		\draw (0.1,1) node [left] {$1$};
		\draw[domain=0.43:5,smooth]
		plot (\x,{0.5*(\x+1)/(\x^2)+1});
		\filldraw [fill=gray!20](0,4.87)--plot [domain=0.43:5,color=white](\x,{0.5*(\x+1)/(\x^2)+1})--(5,1)--(0,1);
		\draw[domain=0:5,color=gray!20] plot (\x,1);
		\draw[domain=0:5,style=dashed] plot (\x,1);
		\draw[domain=0:0.43,color=gray!20] plot (\x,4.87);
		\draw[domain=1:1.13,color=gray!20] plot (5,\x);
		\draw[domain=0:5,smooth] plot (0,\x);
		\end{tikzpicture}
		\caption{The admissible area for $n$ and $p_r-p_1$}
	\end{subfigure}
\caption{}
\end{figure}
{Subsequently,} by using the results above, we can derive the following gradient estimates and Liouville theorems.

\begin{theorem}
	Let $(M^n,g)$ be a complete Riemannian $n$-manifold with Ricci curvature bounded from below by $\operatorname{Ric} \geqslant -K$ {where} $K\geqslant0$, and let $u$ be a positive solution of
	\begin{align}\label{Dpr=auq}
	\tilde{\Delta}_{p_1,...p_r}u + au^q = 0
	\end{align}
	on the ball $B(o,2R)\subset M$ where $1<p_1<...<p_r$ and { $n<\mathcal{N}_1$}. If 
	\begin{align}\label{auq1}
	a>0 \text{ and }\frac{q}{p_1-1}<\frac{n+1}{n-1}+2\sqrt{\frac{1}{(n-1)^2}-\frac{(p_r-p_1)^2}{2(n-1)(p_1-1)^2}},
	\end{align}
	or
	\begin{align}\label{auq2}
	a<0 \text{ and }\frac{q}{p_r-1}>\frac{n+1}{n-1}-2\sqrt{\frac{(p_1-1)^2}{(n-1)^2(p_r-1)^2}-\frac{(p_r-p_1)^2}{2(n-1)(p_r-1)^2}},
	\end{align}
	then  there exists a constant $C$  depending only on $n$, $p_1$, $p_r$, such that
	\begin{align*}
	\frac{|\nabla u|}{u}\leqslant C\frac{1+\sqrt{K}R}{R}
	\end{align*}
	on $B(o,R)$.
	
	In particular, if $M$ is non-compact Riemannian manifold with non-negative Ricci curvature, there is no such positive bounded solution that satisfies (\ref{auq1}) or (\ref{auq2}).
\end{theorem}
\begin{proof}
	Let $\psi(t) = at^{(q-1)/2}$, then $\delta_{\psi}(t) \equiv (q-1)$, and we have known that $\gamma_{\varphi} = \frac{(p_1-1)^2}{n-1}-\frac{(p_r-p_1)^2}{2}$ in Example 4.2. 
	
	Note that either $I_{\psi} = (0,+\infty)$ or $I_\psi = \emptyset$. The former case implies that 
	\begin{align}\label{pr1}
	\frac{n+1}{n-1}-\frac{q}{p_1-1}\geqslant 0 \text{ when } a\geqslant  0,
	\end{align}
	or 
	\begin{align}\label{pr2}
	\frac{n+1}{n-1}-\frac{q}{p_r-1}\leqslant 0\text{ when } a\leqslant  0.
	\end{align}
	The latter one holds if and only if
	$$\sup_{t\geqslant 0}\left(\frac{n+1}{n-1}(\delta_{\varphi}(t)+1)-q\right)^2<\frac{4(p_1-1)^2}{(n-1)^2}-\frac{2(p_r-p_1)^2}{(n-1)},$$
	which infers that 
	\begin{align}\label{pr3_1}
	\frac{n+1}{n-1}(p_r+1)-q<\sqrt{\frac{4(p_1-1)^2}{(n-1)^2}-\frac{2(p_r-p_1)^2}{(n-1)}},
	\end{align}
	and
	\begin{align}\label{pr3_2}
	\frac{n+1}{n-1}(p_1+1)-q>-\sqrt{\frac{4(p_1-1)^2}{(n-1)^2}-\frac{2(p_r-p_1)^2}{(n-1)}}.
	\end{align}
	Combining (\ref{pr1}), (\ref{pr2}), (\ref{pr3_1}) and (\ref{pr3_2}), we obtain the desired results. 
\end{proof}
\begin{remark}
	When equation (\ref{Dpr=auq}) reduces to $p$-Laplacian, then $p_r = p_1 = p$, (\ref{auq1}) and (\ref{auq2}) become the same results in \cite{he2023gradient}.
\end{remark}

It is more interesting to consider what will happen if $I_\psi$ is non-trivial, {in which case, the second critical dimension will make a difference}. Next result shows how these coefficients of equation can affect the set $I_\psi$.
\begin{theorem}
	Let $(M^n,g)$ be a complete Riemannian $n$-manifold with Ricci curvature bounded from below by $\operatorname{Ric} \geqslant -K$ {for some} $K\geqslant0$, and let $u$ be a positive solution of
	\begin{align}\label{Dpr=umus}
	\tilde{\Delta}_{p_1,...p_r}u + u^m-u^k = 0
	\end{align}
	on the ball $B(o,2R)\subset M$ where $1<p_1<...<p_r$ and { $n<\mathcal{N}_1$}. If $m<k$ and
	\begin{align}\label{umus1}
	k\geqslant \frac{n+1}{n-1}(p_r-1)\text{ ~~and~~ }m\leqslant \frac{n+1}{n-1}(p_1-1),
	\end{align}
	then  there exists a constant $C$  depending only on $n$, $p_1$, $p_r$,  such that
	\begin{align*}
	\frac{|\nabla u|}{u}\leqslant C\frac{1+\sqrt{K}R}{R}
	\end{align*}
	on $B(o,R)$.
	
	Furthermore, if {$n<\mathcal{N}_2$}
	then (\ref{umus1}) can be weakened to
	\begin{align}\label{umus2}
	\frac{k}{p_r-1}>\frac{n+1}{n-1}-2\sqrt{\frac{(p_1-1)^2}{(n-1)^2(p_r-1)^2}-\frac{(p_r-p_1)^2}{2(n-1)(p_r-1)^2}},
	\end{align}
	and
	\begin{align}\label{umus3}
	\frac{m}{p_1-1}<\frac{n+1}{n-1}+2\sqrt{\frac{1}{(n-1)^2}-\frac{(p_r-p_1)^2}{2(n-1)(p_1-1)^2}}.
	\end{align}
	
	In particular, if $M$ is non-compact Riemannian manifold with non-negative Ricci curvature, and $u$ is bounded solution, then $u\equiv 1$.
\end{theorem}
\begin{proof}
	Let $\psi(t) = t^{(m-1)/2} - t^{(k-1)/2}$, then $$\delta_{\psi}(t) = \frac{(m-1)t^{(m-1)/2} - (k-1) t^{(k-1)/2}}{t^{(m-1)/2} - t^{(k-1)/2}},$$ and  $d_\varphi = p_r-1$, $l_\varphi = p_1-1$, $\gamma_{\varphi} = \frac{(p_1-1)^2}{n-1}-(p_r-p_1)^2$. 
	Note that 
	\begin{align*}
	I_{\psi}&:= \left\{t>0: \psi(t)\left[\frac{2\left(\delta_\varphi(s) +1 \right)}{n-1} + \delta_\varphi(s)-\delta_\psi(t) \right] \geqslant 0, \text{ for each }s\geqslant 0\right\}\\
	&=\{1\}\cup\left\{t>1:\delta_\psi(t)\geqslant \frac{n+1}{n-1}(p_r-1)-1\right\}\cup\left\{0<t<1:\delta_\psi(t)\leqslant \frac{n+1}{n-1}(p_1-1)-1\right\}.
	\end{align*}
	By De Morgan's laws, we see
	\begin{align*}
	\mathbb{R}^+ -  I_{\psi} = \left\{t>1:\delta_\psi(t)< \frac{n+1}{n-1}(p_r-1)-1\right\}\cup\left\{0<t<1:\delta_\psi(t)> \frac{n+1}{n-1}(p_1-1)-1\right\}.
	\end{align*}

	We then discuss in the following four cases.
	
	\textbf{Case 1:} When $k\geqslant \frac{n+1}{n-1}(p_r-1)$ and $m\leqslant \frac{n+1}{n-1}(p_1-1)$, (\ref{psi2}) naturally holds since $\mathbb{R}^+ - I_{\psi}= \emptyset$.
	
	\textbf{Case 2:} When $\frac{n+1}{n-1}(p_1-1)\leqslant k < \frac{n+1}{n-1}(p_r-1)$ or $ \frac{n+1}{n-1}(p_1-1) < m \leqslant \frac{n+1}{n-1}(p_r-1)$, then from (\ref{psi2}), we have
	\begin{align*}
	\sup_{\substack{s\geqslant 0,\\t\in \mathbb{R}^+ -  I_{\psi}}} \left(\frac{2\left(\delta_{\varphi}(s)+1\right) }{n-1} +\delta_{\varphi}(s) - \delta_{\psi}(t) \right)^2 &= \left(\frac{n+1}{n-1}(p_r-1) - \frac{n+1}{n-1}(p_1-1) \right)^2 \\&< \frac{4(p_1-1)^2}{(n-1)^2}-\frac{2(p_r-p_1)^2}{(n-1)}.
	\end{align*}
	Thus,
	\begin{align}\label{umusc1}
	\frac{(n+1)^2}{2} + (n-1)<\frac{2(p_1-1)^2}{(p_r-p_1)^2} {= \mathcal{N}_1 -1}.
	\end{align}
	{so that $n<\sqrt{2\mathcal{N}_1+3}-2 = \mathcal{N}_2$.}
	
	\textbf{Case 3:} When $m<s<\frac{n+1}{n-1}(p_1-1)$, since
	\begin{align*}
	\sup_{\substack{s\geqslant 0,\\t\in \mathbb{R}^+ -  I_{\psi}}} \left(\frac{2\left(\delta_{\varphi}(s)+1\right) }{n-1} +\delta_{\varphi}(s) - \delta_{\psi}(t) \right)^2 &= \left(\frac{n+1}{n-1}(p_r-1) - k \right)^2 \\&< \frac{4(p_1-1)^2}{(n-1)^2}-\frac{2(p_r-p_1)^2}{(n-1)},
	\end{align*}
	it follows that 
	\begin{align}\label{umusc2}
	\frac{k}{p_r-1}>\frac{n+1}{n-1}-2\sqrt{\frac{(p_1-1)^2}{(n-1)^2(p_r-1)^2}-\frac{(p_r-p_1)^2}{2(n-1)(p_r-1)^2}}.
	\end{align}
	{Note that $k<\frac{n+1}{n-1}(p_1-1)$, hence (\ref{umusc2}) also implies that $n<\mathcal{N}_2$.}
	
	\textbf{Case 4:} When $s>m>\frac{n+1}{n-1}(p_r-1)$, we have
	\begin{align*}
	\sup_{\substack{s\geqslant 0,\\t\in \mathbb{R}^+ -  I_{\psi}}} \left(\frac{2\left(\delta_{\varphi}(s)+1\right) }{n-1} +\delta_{\varphi}(s) - \delta_{\psi}(t) \right)^2 &= \left(\frac{n+1}{n-1}(p_1-1) - m \right)^2 \\&< \frac{4(p_1-1)^2}{(n-1)^2}-\frac{2(p_r-p_1)^2}{(n-1)},
	\end{align*}
	which implies {$n<\mathcal{N}_2$ and}
	\begin{align}\label{umusc3}
	\frac{m}{p_1-1}<\frac{n+1}{n-1}+2\sqrt{\frac{1}{(n-1)^2}-\frac{(p_r-p_1)^2}{2(n-1)(p_1-1)^2}}.
	\end{align}
	
	Combining (\ref{umusc1}), (\ref{umusc2}) and (\ref{umusc3}), we obtain the statements. 
\end{proof}
\begin{remark}\label{pic_umus}
	When $p_r = p_1 = p = 2$, {so that $\mathcal{N}_2 = \infty$ and $n <  \mathcal{N}_2$ naturally} holds. Thus the (\ref{umus2}) and (\ref{umus3}) show that
	$$m<\frac{n+3}{n-1} \text{ ~~and~~ }k>1,$$
	which improve Wang's result in \cite{wang2024gradient} (see Figure 3):
	$$1<m<\frac{n+3}{n-1}\text{ ~~or~~ }1<k<\frac{n+3}{n-1}.$$
\end{remark}
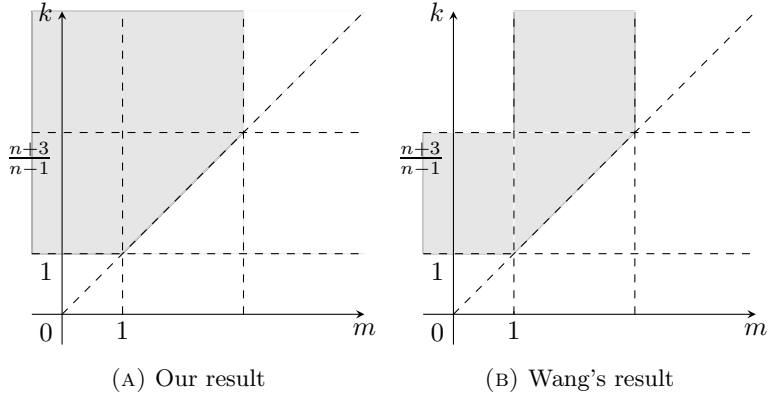
\begin{figure}[H]
	\begin{subfigure}{0.4\textwidth}
		\begin{tikzpicture}[scale=0.8]
		\filldraw [fill=gray!20] (-0.5,1) -- (1,1) -- (3,3) -- (3,5) -- (-0.5,5) -- (-0.5,1);
		\draw [domain=-0.5:5,smooth,color=gray!20] plot (\x,5);
		\draw [domain=1:5,smooth,color=gray!20] plot (-0.5,\x);
		\draw [domain=-0.5:5,smooth,color=gray!20] plot (\x,1);
		\draw [domain=1:5,smooth,color=gray!20] plot (\x,\x);
		\draw [domain=1:5,smooth,color=gray!20] plot (3,\x);
		%%%%
		\draw[-stealth] (-0.5,0)--(5,0) node[below]{$m$};
		\draw[-stealth] (0,-0.5)--(0,5) node[left]{$k$};
		\draw (0,0) node [below left] {$0$};
		\foreach \i in {1}{\draw (\i,0)--node [below] {$1$}(\i,0.05);}
		\draw (0,1) node [below left] {$1$};
		\draw (0,3) node [below left] {$\frac{n+3}{n-1}$};
		\draw [domain=3:5,smooth,color=white] plot   (\x,5);
		\draw[domain=-0.5:5,style=dashed] plot (\x,1);
		\draw[domain=-0.5:5,style=dashed] plot (\x,3);
		\draw[domain=0:5,style=dashed,variable=\t] plot (1,\t);
		\draw[domain=0:5,style=dashed,variable=\t] plot (3,\t);
		\draw[domain=0:5,style=dashed,variable=\t] plot (\t,\t);
		\end{tikzpicture} 
		\caption{Our result}
	\end{subfigure}
	\begin{subfigure}{0.4\textwidth}
		\begin{tikzpicture}[scale=0.8]
		\filldraw [fill=gray!20] (-0.5,1) -- (1,1) -- (3,3) -- (3,5) -- (1,5) -- (1,3) -- (-0.5,3) -- (-0.5,1);
		\draw [domain=-0.5:5,smooth,color=white] plot (\x,5);
		\draw [domain=-0.5:5,smooth,color=white] plot (-0.5,\x);
		\draw [domain=-0.5:5,smooth,color=white] plot (\x,1);
		\draw [domain=-0.5:5,smooth,color=white] plot (\x,\x);
		\draw [domain=-0.5:5,color=white] plot (3,\x);
		\draw [domain=-0.5:5,color=white] plot (1,\x);
		\draw [domain=-0.5:5,color=white] plot (\x,3);
		%%%%
		\draw[-stealth] (-0.5,0)--(5,0) node[below]{$m$};
		\draw[-stealth] (0,-0.5)--(0,5) node[left]{$k$};
		\draw (0,0) node [below left] {$0$};
		\foreach \i in {1}{\draw (\i,0)--node [below] {$1$}(\i,0.05);}
		\draw (0,1) node [below left] {$1$};
		\draw (0,3) node [below left] {$\frac{n+3}{n-1}$};
		\draw [domain=3:5,smooth,color=white] plot   (\x,5);
		\draw[domain=-0.5:5,style=dashed] plot (\x,1);
		\draw[domain=-0.5:5,style=dashed] plot (\x,3);
		\draw[domain=0:5,style=dashed,variable=\t] plot (1,\t);
		\draw[domain=0:5,style=dashed,variable=\t] plot (3,\t);
		\draw[domain=0:5,style=dashed,variable=\t] plot (\t,\t);
		\end{tikzpicture}
		\caption{Wang's result}
	\end{subfigure}
	\caption{The admissible areas for Liouville theorem compared with \cite{wang2024gradient}}
\end{figure}
\begin{remark}
	Shortly after the completion of this manuscript, we saw a new paper \cite{he2024localgloballoggradientestimates} uploaded on arXiv by J. He and his collaborators, which provides a gradient estimate of the equation
	$$\Delta_pu+bu^q+cu^r =0,$$
	This is also a special case of Theorem \ref{mainthm}, by taking $\varphi(t)= t^{p-2}$ and $\psi(t) = bt^{\frac{q-1}{2}} + ct^{\frac{r-1}{2}}$.
\end{remark}

One might ask what if $\psi$ is not a polynomial. To illustrate this, we will give the following theorem.
\begin{theorem}\label{log}
	Let $(M^n,g)$ be a complete Riemannian $n$-manifold with Ricci curvature bounded from below by $\operatorname{Ric} \geqslant -K$ {where} $K\geqslant0$, and let $u$ be a positive solution of
	\begin{align}\label{Dpr=log}
	\tilde{\Delta}_{p_1,...p_r}u + au^q(\log u)^m = 0
	\end{align}
	on the ball $B(o,2R)\subset M$ where $m = \frac{2k_1+1}{2k_2+1}$ where $k_1$ and $k_2$ are integers, $ma<0$, $1<p_1<...<p_r$, and {$n<\mathcal{N}_2$}.
	If
	\begin{align}
	\frac{q}{p_r-1}>\frac{n+1}{n-1}-2\sqrt{\frac{(p_1-1)^2}{(n-1)^2(p_r-1)^2}-\frac{(p_r-p_1)^2}{2(n-1)(p_r-1)^2}},
	\end{align}
	and
	\begin{align}
	\frac{q}{p_1-1}<\frac{n+1}{n-1}+2\sqrt{\frac{1}{(n-1)^2}-\frac{(p_r-p_1)^2}{2(n-1)(p_1-1)^2}},
	\end{align}
	then  there exists a constant $C$  depending only on $n$, $p_1$, $p_r$, such that
	\begin{align*}
	\frac{|\nabla u|}{u}\leqslant C\frac{1+\sqrt{K}R}{R}~~~~~
	\end{align*}
	on  $B(o,R)$.	
	
	In particular, if $M$ is non-compact Riemannian manifold with non-negative Ricci curvature, and $u$ is bounded solution, then $u\equiv 1$ when $m>0$, there is no such positive solution when $m<0$.
\end{theorem}
\begin{proof}
	Let $\psi(t) = at^{(q-1)/2}\left(\frac{1}{2}\log t\right)^m$, then $$\delta_{\psi}(t) = (q-1) + \frac{2m}{\log t}.$$ 
	When $a>0$, similarly
	\begin{align*}
	I_{\psi}
	=\{1\}\cup\left\{t>1:\delta_\psi(t)\geqslant \frac{n+1}{n-1}(p_r-1)-1\right\}\cup\left\{0<t<1:\delta_\psi(t)\leqslant \frac{n+1}{n-1}(p_1-1)-1\right\},
	\end{align*}
	and
	\begin{align*}
	\mathbb{R}^+ -  I_{\psi} = \left\{t>1:\delta_\psi(t)< \frac{n+1}{n-1}(p_r-1)-1\right\}\cup\left\{0<t<1:\delta_\psi(t)> \frac{n+1}{n-1}(p_1-1)-1\right\}.
	\end{align*}
	
	Now, we discuss in the following three cases.
	
	\textbf{Case 1:} When $q\geqslant \frac{n+1}{n-1}(p_r-1)$, we see
	\begin{align*}
	\sup_{\substack{s\geqslant 0,\\t\in \mathbb{R}^+ -  I_{\psi}}} \left(\frac{2\left(\delta_{\varphi}(s)+1\right) }{n-1} +\delta_{\varphi}(s) - \delta_{\psi}(t) \right)^2 &= \left(q - \frac{n+1}{n-1}(p_1-1) \right)^2 \\&< \frac{4(p_1-1)^2}{(n-1)^2}-\frac{2(p_r-p_1)^2}{(n-1)},
	\end{align*}
	which implies
	\begin{align}\label{log_c1}
	\frac{q}{p_1-1}<\frac{n+1}{n-1}+2\sqrt{\frac{1}{(n-1)^2}-\frac{(p_r-p_1)^2}{2(n-1)(p_1-1)^2}}.
	\end{align}
	{Since $q\geqslant \frac{n+1}{n-1}(p_r-1)$, it must hold that
		\begin{align}\label{log_c2}
		\frac{(n+1)^2}{2} + (n-1)<\frac{2(p_1-1)^2}{(p_r-p_1)^2},
		\end{align}	
thus $n<\mathcal{N}_2$.} 
	
	\textbf{Case 2:} When $\frac{n+1}{n-1}(p_1-1) <q < \frac{n+1}{n-1}(p_r-1)$, then from (\ref{psi2}), we have
	\begin{align*}
	\sup_{\substack{s\geqslant 0,\\t\in \mathbb{R}^+ -  I_{\psi}}} \left(\frac{2\left(\delta_{\varphi}(s)+1\right) }{n-1} +\delta_{\varphi}(s) - \delta_{\psi}(t) \right)^2 &= \left(\frac{n+1}{n-1}(p_r-1) - \frac{n+1}{n-1}(p_1-1) \right)^2 \\&< \frac{4(p_1-1)^2}{(n-1)^2}-\frac{2(p_r-p_1)^2}{(n-1)}.
	\end{align*}
	Hence,
	\begin{align*}
	\frac{(n+1)^2}{2} + (n-1)<\frac{2(p_1-1)^2}{(p_r-p_1)^2},
	\end{align*}
	{and then $n<\mathcal{N}_2$.}
	
	\textbf{Case 3:} When $q\leqslant \frac{n+1}{n-1}(p_1-1)$, since
	\begin{align*}
	\sup_{\substack{s\geqslant 0,\\t\in \mathbb{R}^+ -  I_{\psi}}} \left(\frac{2\left(\delta_{\varphi}(s)+1\right) }{n-1} +\delta_{\varphi}(s) - \delta_{\psi}(t) \right)^2 &= \left(\frac{n+1}{n-1}(p_r-1) - q \right)^2 \\&< \frac{4(p_1-1)^2}{(n-1)^2}-\frac{2(p_r-p_1)^2}{(n-1)},
	\end{align*}
	it follows that {$n<\mathcal{N}_2$ and}
	\begin{align}\label{log_c3}
	\frac{q}{p_r-1}>\frac{n+1}{n-1}-2\sqrt{\frac{(p_1-1)^2}{(n-1)^2(p_r-1)^2}-\frac{(p_r-p_1)^2}{2(n-1)(p_r-1)^2}}.
	\end{align}
	
	Combining (\ref{log_c1}), (\ref{log_c2}) and (\ref{log_c3}), we finish the proof.
\end{proof}
\begin{remark}
	When $\varphi \equiv 1$, B. Peng \cite{peng2021yau} gave a gradient estimate for
	$a\neq 0$ and $m = \frac{k_1}{2k_2+1}>2$, although
	Theorem \ref{log} requires $k_1$ to be odd, our result is still feasible for $m<2$, even $m$ is negative. Moreover, the gradient estimate in \cite{peng2021yau} is not Cheng--Yau-type, which cannot derive the Liouville property of {that} equation.
\end{remark}

\section*{Acknowledge}
The author is very grateful to the reviewers for their careful review and valuable comments.

%\bibliography{refsab}
%\bibliographystyle{abbrv}

\end{document}